\pdfoutput=1

\documentclass[a4paper,12pt]{amsart}
\newcommand{\topcounter}{section}

\usepackage[applemac]{inputenc}

\usepackage[a4paper]{geometry}
\usepackage{graphicx}
\usepackage{amsmath}
\usepackage{amsthm}
\usepackage{amssymb}
\usepackage{amsfonts}
\usepackage{mathrsfs}
\usepackage{enumerate}

\usepackage{caption}
\theoremstyle{plain}
\newtheorem{theorem}{Theorem}[\topcounter]

\newtheorem{dictionnary*}{Dictionary}

\newtheorem{lemma}[theorem]{Lemma}
\newtheorem*{lemma*}{Lemma}

\newtheorem*{corollary*}{Corollary}
\newtheorem*{important}{Important}

\theoremstyle{definition}
\newtheorem{definition}{Definition}[\topcounter]

\newtheorem*{demonstration}{Proof}

\newtheorem*{convention*}{Convention}

\newtheorem*{generalisation*}{Generalisation}

\newtheorem*{application*}{Application}

\newtheorem*{remark}{Remark}

\newtheorem*{variation*}{Variation}

\def\inv{^{-1}}

\DeclareMathOperator{\id}{id}

\newcommand{\confer}{\textit{cf.} }

\newcommand{\ZZ}{\mathbb{Z}}

\newcommand{\Map}{\Lambda}
\newcommand{\Gr}{\Gamma}

\newcommand{\shortpar}{\vspace{-5mm}}

\renewcommand{\epsilon}{\varepsilon}

\usepackage[algoruled, vlined, algosection]{algorithm2e}
\SetKwComment{comment}{$(*$\:}{\:$*)$}

\renewcommand{\phi}{\varphi}

\newcommand{\fra}{\mathfrak{a}}
\newcommand{\frb}{\mathfrak{b}}
\newcommand{\frc}{\mathfrak{c}}
\newcommand{\frd}{\mathfrak{d}}
\newcommand{\fre}{\mathfrak{e}}
\newcommand{\frf}{\mathfrak{f}}

\usepackage[all,cmtip]{xy}


\begin{document}

\title[An Optimal Algorithm to Generate rooted Trivalent Diagrams and rooted Triangular Maps]{An Optimal Algorithm to Generate rooted Trivalent Diagrams and rooted Triangular Maps}
\author{Samuel Alexandre Vidal\\Lille 1 University of Sciences and Technology\\Paul Painlev Laboratory of Mathematics\\59 655 Villeneuve d'Ascq Cdex France}

\maketitle

\begin{abstract}
A \emph{trivalent diagram} is a connected, two-colored bipartite graph (parallel edges allowed but not loops) such that every black vertex is of degree $1$ or $3$ and every white vertex is of degree $1$ or $2$, with a cyclic order imposed on every set of edges incident to to a same vertex.
A \emph{rooted} trivalent diagram is a trivalent diagram with a distinguished edge, its \emph{root}.
We shall describe and analyze  an algorithm giving an exhaustive list of rooted trivalent diagrams of a given size (number of edges), the list being non-redundant in that no two diagrams of the list are isomorphic.
The algorithm will be shown to have optimal performance in that the time necessary to generate a diagram will be seen to be bounded in the amortized sense, the bound being independent of the size of the diagrams. That's what we call the CAT property.
One objective of the paper is to provide a reusable theoretical framework for algorithms generating exhaustive lists of complex combinatorial structures with attention paid to the case of unlabeled structures and to those generators having the CAT property.
\end{abstract}

\setcounter{section}{-1}

\parskip = 0 mm
\parindent = 5 mm
\sloppy

\tableofcontents

\parskip = 5 mm
\parindent = 5 mm
\sloppy

\section{Introduction}

Roughly speaking, a trivalent diagram is a connected graph with degree conditions on its vertices and cyclic orientations on the edges adjacent to each vertex
It is the combinatorial description of an unembedded trivalent ribbon graph \cite{thurston82,kontsevich92} (\confer definitions \ref{def:triv:diag} and \ref{def:morp:triv:diag} for a precise definition). We shall see (\confer theorem \ref{th:reconstruction:triv:diag:perm}) that it can be described by a pair of permutations $(\sigma_\bullet,\sigma_\circ)$ satisfying the conditions of \emph{involutivity} $\sigma_\circ^2 = \id$ and \emph{triangularity} $\sigma_\bullet^3 = \id$. The notion of rooted trivalent diagrams is also very useful, both to our study and to the target applications; so we take a special care to study them in detail.

\subsection{Problem Statement}

We shall describe and analyze two algorithms, the first giving an exhaustive list of rooted trivalent diagrams of a given size (\confer definitions \ref{def:root:triv:diag} and \ref{def:morp:root:triv:diag} below) and the second giving an exhaustive list of unrooted trivalent diagrams (definitions \ref{def:triv:diag} and \ref{def:morp:triv:diag}), those lists being non-redundant in that no two diagrams of the lists are isomorphic.
The algorithm for rooted diagrams will be shown to have optimal performance meaning that the time necessary to generate a diagram is bounded in the amortized sense. What is striking is that the bound is independent of the size of the diagrams being generated.
One objective of the paper is to provide a reusable theoretical framework for algorithms generating exhaustive lists of complex combinatorial structures with attention paid to the case of unlabeled structures and to those generators having the CAT property.

\subsection{Motivations}

In a recent paper \cite{vidal06} we gave a complete classification of the subgroups of the modular group $\mathrm{PSL}_2(\ZZ)$ and their conjugacy classes by rooted trivalent diagrams and trivalent diagrams. A question one may ask is how to generate a complete list of such trivalent diagrams. Such a question is unavoidable: for a classification to be fully satisfactory there should be a systematic way to enumerate all the particular instances of the objects being classified.
Moreover, it was soon realized that there is a connection with combinatorial maps. In this paper we clarify that point and give as an application a way to generate exhaustive lists of triangular combinatorial maps.

The other sources of motivations to generate trivalent diagrams come mainly from mathematical physics in connection with two-dimensionnal quantum gravity and the Witten-Kontsevich model \cite{kontsevich92}. Algebraic topology is also a source of motivation through triangular subdivisions of surfaces, knots, braids, links and tangles theory \cite{thurston82, barnathan98}. It is also connected to the deformation theory of quantized Hopf algebras \cite{drinfeld90, drinfeld91}. The problem we solve is also relevant to the study of combinatorial maps as explained in section \ref{sect:application:triangular:maps} and to the vast galoisian program of A. Grothendieck \cite{grothendieckprogramme} as explained in hundreds of papers and books such as \cite{schneps94, vidal06, landozvonkine04, douadydouady}. As an application, we give in section \ref{sect:application:modular:group} a way to generate a complete list describing all the sub-groups of a given finite index in the modular group $\mathrm{PSL}_2(\ZZ)$ and a way to decide conjugacy relations among those subgroups. We show also in section \ref{sect:application:triangular:maps}, as a second application, how to generate an exhaustive list of triangular maps satisfying various criteria.

\subsection{What is a CAT Generator}

The expression ``CAT generator" is an acronym for \emph{constant amortized time generator} meaning a generator of combinatorial structures that average spend only a constant time generating each of the structures.
The usual idea in such a generator is that passing from one structure to the next requires only a few modifications to be made. Sometimes, though, it could take more modifications than usual and we don't usually have any upper bound on the amount of actual modifications that could be needed to pass from a structure to the next. When need for large number of modifications tends to be significantly rare in comparison to small ones, we can sometimes prove that an amortization effect is going on.
Technically, one can summarize that amortization effect in saying that the total amount of time needed to generate $n$ distinct structures is asymptotically bounded by a constant multiple of the number $n$ of structures being generated, the word \emph{constant} meaning that the bound is independent of the size of the structures being generated.


\section{Trivalent Diagrams}

\begin{definition}
\label{def:triv:diag}
A trivalent diagram is a connected, two colored bipartite graph (parallel edges allowed but not loops) such that every black vertex is of degree $1$ or $3$ and every white vertex is of degree $1$ or $2$, with a cyclic order imposed on the edges incident to each vertex. The \emph{size} of a trivalent diagram is the number of its edges.
\end{definition}

Given a trivalent diagram $\Gr$, we denote by $\Gr_{-}$, $\Gr_\bullet$ and $\Gr_\circ$ the sets of its edges, black vertices and white vertices respectively. Given an edge $a \in \Gr_{-}$, we denote by $\partial_\bullet(a)\in \Gr_\bullet $ and $\partial_\circ(a)\in \Gr_\circ$ the black vertex and the white vertex to which it is incident.
Given an edge $a \in \Gr_{-}$, we denote by $\sigma_\bullet(a)$ and $\sigma_\circ(a)$ the next edge incident to $\partial_\bullet(a)$ and $\partial_\circ(a)$ respectively, in the cyclic order. According to the degree conditions of the definition we have $\sigma_\bullet^3=\sigma_\circ^2 = \id$ which implies that both $\sigma_\bullet$ and $\sigma_\circ$ are bijections, so they are permutations on the set of edges $\Gr_{-}$ of the diagram. The connectivity condition of the definition is equivalent to the transitivity of the permutation group generated by $\sigma_\bullet$ and $\sigma_\circ$.

\begin{figure}
\begin{center}
\includegraphics{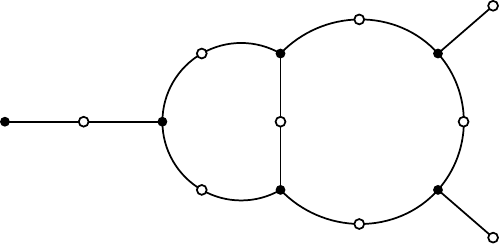}
\caption{A trivalent diagram is conveniently described by a little diagram like the one above, hence the name. The actual cyclic orientation of the vertices are conveniently rendered implicit by adopting the counter-clockwise orientation of the figure.}
\end{center}
\end{figure}

\begin{definition}
\label{def:morp:triv:diag}
A \emph{morphism} $\phi$ between two trivalent diagrams $\Gr$ and $\Gr'$ is a triple of mappings
$\phi_\bullet$, $\phi_\circ$ and $\phi_{-}$ from the three sets 
$\Gr_\bullet$, $\Gr_\circ$ and $\Gr_{-}$ to the three sets $\Gr'_\bullet$, $\Gr'_\circ$ and $\Gr'_{-}$ respectively,
compatible with the three structure mappings and the two permutaion in that $\phi_{-}$ is both equivarient to the $\sigma_\bullet$ and $\sigma_\circ$ permutations and the following diagram is commutative.
\begin{align*}
\xymatrix@C=1.5cm@R=1.5cm
{
	{\Gr_\bullet}
		\ar[d]_{\phi_\bullet}			&
	{\Gr_{-}}
		\ar[l]_{\partial_\bullet}
		\ar[r]^{\partial_\circ}
		\ar[d]^{\phi_{-}}				&
	{\Gr_\circ}
		\ar[d]^{\phi_\circ}			\\
	{\Gr'_\bullet}					&
	{\Gr'_{-}}
		\ar[l]^{\partial_\bullet}
		\ar[r]_{\partial_\circ}			&
	{\Gr'_\circ}
}
\end{align*}
When those three mappings are bijections the morphism is an \emph{isomorphism}.
\end{definition}

A important fact, well known to the experts is recalled in the following theorem. It is used throughout the article to formulate the algorithms.

\begin{theorem}
\label{th:reconstruction:triv:diag:perm}
The set $\Gr_{-}$ and the two permutions $\sigma_\bullet$ and $\sigma_\circ$ altogether suffice to describe the isomorphism class of the digram $\Gr$. Moreover, the cycles of the permutations $\sigma_\bullet$ and $\sigma_\circ$ are in natural bijection with the black and white vertice of $\Gr$ respectively.
\end{theorem}

\begin{demonstration}
Given a trivalent diagram $\Gr$, an isomorphic trivalent diagram $\Gr'$ can be reconstructed from the set $\Gr_{-}$ and the two permutations $\sigma_\bullet$ and $\sigma_\circ$ of $\Gr_{-}$. Let $\Gr_{-}/\sigma_\bullet$, the cyles of the permutation $\sigma_\bullet$, be its set of black vertice and let $\Gr_{-}/\sigma_\circ$, the cycles of the permutation $\sigma_\circ$, be its set of white vertice, and define its boundary mappings $\partial'_\bullet$ and $\partial'_\circ$ to be the natural projection of the quotients.

We now construct the isomorphism the following way. Since $\partial_\bullet$ and $\partial_\circ$ are equivarient to the permutations $\sigma_\bullet$ and $\sigma_\circ$ respectively, they induce natural maps $\phi_\bullet$ and $\phi_\circ$ completing the following commutative diagram.
\begin{align*}
\xymatrix@C=1.5cm@R=1.5cm
{
	{\Gr'_\bullet}
		\ar@{.>}[d]_{\phi_\bullet}			&
	{\Gr'_{-}}
		\ar[l]_{\partial'_\bullet}
		\ar[r]^{\partial'_\circ}
		\ar@{=}[d]^{\phi_{-}}				&
	{\Gr'_\circ}
		\ar@{.>}[d]^{\phi_\circ}			\\
	{\Gr_\bullet}						&
	{\Gr_{-}}
		\ar[l]^{\partial_\bullet}
		\ar[r]_{\partial_\circ}				&
	{\Gr_\circ}
}
\end{align*}
Taking $\phi_{-}$ to be the identity mapping, one has a morphism from the diagram $\Gr'$ to the diagram 
$\Gr$. To show it's an isomorphism one has to show the bijectivity of the three mappings $\phi_\bullet$, $\phi_\circ$ and $\phi_{-}$. The mapping $\phi_{-}$, being the identity, is necessarily bijective. The bijectivity of the mappings $\phi_\bullet$ and $\phi_\circ$ means that two edges are in the same cycles of the respective permutations $\sigma_\bullet$ and $\sigma_\circ$ if and only if they are incident to the same black and white vertex respectively, which is guaranteed by the definition.
\hfill$\Box$
\end{demonstration}

\subsection{Rooted Trivalent Diagrams}

The following notion play an important rle in that dissertation and in the applications.

\begin{definition}
\label{def:root:triv:diag}
A trivalent diagram is said to be \emph{rooted} if one of its edges is distinguished from the others as its \emph{root}.
\end{definition}
\shortpar

A convenient way to describe the rooting of a diagram is to draw a little cross on its distinguished edge.

\begin{definition}
\label{def:morp:root:triv:diag}
A \emph{morphism} $\phi$ of rooted trivalent diagrams $(\Gr,a)$ and $(\Gr',a')$ is a morphism of the underlying diagrams (ignoring the roots) which $\phi_{-}$ component is further assumed to send root to root.
\end{definition}

\subsection{Labeled \textit{versus} Unlabeled Diagrams}

Historically, the dichotomy between \emph{labeled} and \emph{unlabeled} structures had been greatly clarified and properly emphasized by the introduction by A. Joyal of \emph{combinatorial species} \cite{Joyal81}. The subject was, and still is, a very prolific source of discovery from the Quebec school of combinatorics and from a growing community of researchers around the world. One must cite the wonderful book \cite{bergeron94, bergeron98} by  F. Bergeron, G. Labelle and P. Leroux, which gives a briliant exposition of the whole subject.

On a given set of vertices $X$ one can build different trivalent diagrams and rooted trivalent diagrams. We denote by $D_3(X)$ and $D^\bullet_3(X)$ the corresponding sets of structures, we call $X$ the \emph{labeling alphabet} and we refer to diagrams one can build on that set as diagrams \emph{labeled} by $X$.
Any bijection $\varrho$ between two finite sets $X$ and $Y$ induces a bijection $\varrho_*$ between the sets $D_3(X)$ and $D_3(Y)$ of trivalent diagrams labeled by $X$ and $Y$, respectively. This induced bijection is the \emph{relabeling operation} from $D_3(X)$ to $D_3(Y)$. It is also referred as a \emph{transport of structure} along the relabeling bijection $\varrho$. The same considerations also apply to rooted trivalent diagrams and in fact to any labeled combinatorial structures.

The above discution leads to the consideration of the \emph{Joyal Functors} $D_3$ and $D^\bullet_3$ of the two combinatorial species of trivalent diagrams and rooted trivalent diagrams, respectively.
In the formalism of Joyal, two labeled structures are said to be \emph{conjugate} or \emph{isomorphic} if they coincide modulo the relabeling operation.
An unlabeled structure is then just a conjugacy class of labeled structures. We denote by $\tilde{D}_3(n)$ and $\tilde{D}^\bullet_3(n)$ the sets of unlabeled trivalent diagrams, unrooted and rooted respectively, and by $D_3(n)$ and $D^\bullet_3(n)$ the sets of trivalent diagrams, unrooted and rooted respectively and labeled by the set  $\{\,1, \dots, n\,\}$. The symmetric group $\mathfrak{S}_n$ acts via relabeling on the set of structures labeled by $\{\,1, \dots, n\,\}$ and the sets $\tilde{D}_3(n)$ and  $\tilde{D}^\bullet_3(n)$ can be seen as the quotient sets of those group actions.
\begin{align*}
	 \tilde{D}_3(n) \overset{\text{def.}}{=}  D_3(n) / \mathfrak{S}_n
	 \quad \text{ and } \quad
	 \tilde{D}^\bullet_3(n) \overset{\text{def.}}{=}  D^\bullet_3(n) / \mathfrak{S}_n
\end{align*}
We call the corresponding natural projections,
\begin{align*}
	\pi_n &:  D_3(n) \to \tilde{D}_3(n) \\
	\pi_n &: D^\bullet_3(n) \to  \tilde{D}^\bullet_3(n)
\end{align*}
the \emph{condensation mappings} of the combinatorial species $D_3$ and $D^\bullet_3$.

\section{Characteristic Labeling}
\label{sect:char:label}

A \emph{characteristic labeling} is the choice of a unique representative in every conjugacy class of structures. In other terms, a characteristic labeling can be seen as a natural section to the condensation mapping $\pi$.
The characteristic labelings that we like are those which are computable. We like them even more if there is an efficient way to compute them.

Rooted trivalent diagrams have the enjoyable property of possessing many characteristic labelings that are computable by means of efficient algorithms. This situation is to be contrasted with that of general graphs. No algorithm is known to decide in polynomial time whether two given graphs are isomorphic, and having an efficient algorithm computing characteristic labeling of general graphs would render that particular problem trivial. What makes trivalent diagrams particular in that respect is not so much that they are trivalent but more in that their edges are cyclically oriented at their vertices. Indeed, general graphs with only trivalent vertices still suffer from the above problem.

What we give now, is a succinct description of an algorithm producing a characteristic labeling of rooted trivalent diagrams $\Gr$ and having linear time complexity in the number of edges of $\Gr$.
The idea is the following: build a rooted planar binary tree $T$ by depth-first traversal (in prefix order) of the edges of the diagram (not the vertices, we insist on the edges). Given a particular edge $a$ of $\Gr$, the two directions that are explored from it, are given by the two operations $\sigma_\bullet$ and $\sigma_\circ$ on the set of edges. We take care never to revisit a previously visited edge and we label the edges of $\Gr$ by numbers from $1$ to $n$ according to the order of their appearance in the depth-first traversal.
\begin{figure}
\begin{center}
\includegraphics{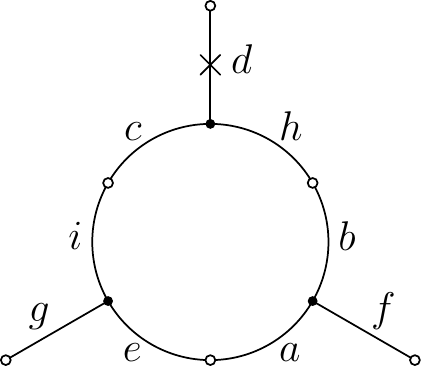}
\hspace{2cm}
\includegraphics{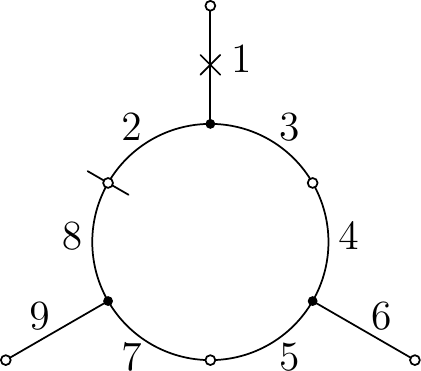}
\caption{If one gives as input to the relabeling algorithm (algorithm \ref{algo:relabel}) the rooted diagram shown on the left with an arbitrary initial labeling on the arbitrary alphabet $X=\{\,a,b,c,d,e,f,g,h,i\,\}$, it produces the characterstic relabeling shown on the right with numbers from $1$ to $9$ according to the depth-first traversal order of algorithm \ref{algo:visit:relabel}. Note the \emph{natural cutting} between the edges labeled by $2$ and $8$ that arises from the depth-first traversal.}
\end{center}
\end{figure}

\subsection{Implementation}
\linesnumbered

We need as global data an integer $c$ and the following seven arrays:
\begin{align*}
	visited &: X \to \mathrm{Bool}		\\
	\ell_0 &: X \to \{\, 1, ..., n \,\}		\\
	\ell_1 &: \{\, 1, ..., n \,\} \to X		\\
	s_0, s_1 &: X \to X				\\
	 t_0, t_1  &: \{\, 1, ..., n \,\} \to \{\, 1, ..., n \,\}.
\end{align*}
Algorithm \ref{algo:visit:relabel} which is an auxiliary recursive program computes the transport bijections $\ell_0$ and $\ell_1$. Algorithm \ref{algo:relabel} is the main entry point of the relabeling process. It does the initialization job (line 2 to 4) and the actual relabeling of the input diagram (line 6 to 8). It takes as input a trivalent diagram labeled with the elements of the set $X$ and rooted by the element $x$ of $X$. The arrays $s_0$ and $s_1$ and the element $x \in X$ are descriptions of the input diagram via its associated two permutations $\sigma_\bullet$ and $\sigma_\circ$ (\confer theorm \ref{th:reconstruction:triv:diag:perm}). The output diagram is encoded by the two arrays $t_0$ and $t_1$ in the very same fashion. The $visited$ array is used to remember the positions already visited by the relabeling process. The integer $c$ serves as a counter to label the vertices in the order they are encountered, $\ell_0$ and $\ell_1$ are internal arrays describing the mutual inverse transport bijections between the input diagram and the output diagram.

\begin{algorithm}
\dontprintsemicolon
\Begin{
	\lIf{$visited \,[\,x\,]$}{return} \;
	$visited \,[\,x\,] \leftarrow true$ \;
	$\ell_0 \,[\,x\,] \leftarrow c$ \;
	$\ell_1 \,[\,c\,] \leftarrow x$ \;
	$c \leftarrow c+1$ \;
	$\textsc{Visit }(s_0 \,[\,x\,] )$ \;
	$\textsc{Visit }(s_1 \,[\,x\,] )$ \;
}
\caption{$\textsc{Visit }(x : X)$}
\label{algo:visit:relabel}
\end{algorithm}
\begin{algorithm}
\dontprintsemicolon
\Begin{
	$c \leftarrow 1$ \;
	\For{$i \in X$}{
		$visited \,[\,i\,] \leftarrow false$ \;	
	}
	$\textsc{Visit } (x)$ \;
	\For{$k \in \{\, 1, ..., n \,\}$}{
		$t_0 \,[\,k\,] \leftarrow \ell_0 \,[\,s_0 \,[\, \ell_1 \,[\,k\,]\,]\,]$ \;	
		$t_1 \,[\,k\,] \leftarrow \ell_0 \,[\,s_1 \,[\, \ell_1 \,[\,k\,]\,]\,]$ \;	
	}
}
\caption{$\textsc{Relabel }(x : X)$}
\label{algo:relabel}
\end{algorithm}

\subsection{Correctness}

The idea behind that algorithm is quite simple and presents no difficulty excepting the actual proof of the relabeling being characteristic. There are two ways to do the proof~; one is conceptual by nature and the other is more technical. The particular description of the algorithm is itself part of that former argument. We shall give both arguments because preferring one or the other is simply a mater of taste. Let's give the conceptual argument first.

One could have taken the input diagram to be labeled by the set $\{\,1,...,n\,\}$ and then shown that the output labeled diagram remains unchanged if one conjugates the input labeled diagram according to any permutation of the labeling set. Such a proof would typically look rather technical if not difficult. Instead, one can rather \emph{abstract} the labeling alphabet of the input diagram to be an arbitrary $n$ elements set $X$, this requirement being the only assumption made on $X$.
In particular, we make absolutely no assumption on its elements or on any structure that it may carry.

A moment's thought may convince the reader that abstracting the input label set to $X$ and making no assumption whatsoever on its elements indeed guarantees the required invariance. As this argument is a bit subtle and may seems a hand-waving argument to most people~; we now give another proof avoiding such considerations.

\begin{theorem}
\label{theorem:characteristic}
Algorithm {\rm \ref{algo:relabel}} produces a characteristic relabeling of the connected rooted trivalent diagrams of size $n$ -- that is, $t_0$ and $t_1$ are invariant under any bijection from $X$ onto another set $X'$.
\end{theorem}

\begin{demonstration}
Any bijection $\varrho$ between two sets of input labels $X$ and $X'$ induces a conjugacy of the two input permutations $s_0$ and $s_1$ of $X$ yielding two permutations $s_0' = \varrho\cdot s_0\cdot\varrho\inv$ and $s_1' = \varrho\cdot s_1\cdot\varrho\inv$ of $X'$.
Now, putting $\ell_0 (x) = c$ and $\ell_0'(x')=c$ according to line 4 of Algorithm \ref{algo:visit:relabel} with $x' = \varrho (x)$ and varying $x$ yields $\ell_0' = \ell_0 \cdot \varrho\inv$. Similarly, considering line 5 of the same algorithm, we get $\ell_1' = \varrho \cdot \ell_1$. The permutations $t_0$, $t_0'$, $t_1$ and $t_1'$ verify the following identities (by line 6-8 of algorithm \ref{algo:relabel}):
\begin{align*}
t_0 &= \ell_0 \cdot s_0 \cdot \ell_1	,	&
t_0' &= \ell_0' \cdot s_0' \cdot \ell_1'	,	\\
t_1 &= \ell_0 \cdot s_1 \cdot \ell_1	,	&
t_1' &= \ell_0' \cdot s_1' \cdot \ell_1',
\end{align*}
and substituting $\ell_0'$, $\ell_1'$, $s_0'$, and $s_1'$ for their above values yields, a cancellation of the $\varrho$'s,
\begin{align*}
	t_0' &= (\ell_0 \cdot \varrho\inv)\cdot(\varrho \cdot s_0 \cdot \varrho\inv)\cdot(\varrho \cdot \ell_1) = t_0,		\\
	t_1' &= (\ell_0 \cdot \varrho\inv)\cdot(\varrho \cdot s_1 \cdot \varrho\inv)\cdot(\varrho \cdot \ell_1) = t_1,
\end{align*}
thus proving the required invariance of the output.\hfill$\Box$
\end{demonstration}




\section{Generating Algorithm}
\label{sect:generator}

Let's imagine that while exploring a particular rooted trivalent diagram using algorithm \ref{algo:relabel} of section \ref{sect:char:label}, we output a sequence of events describing the particular cycles of the permutations $t_0$ and $t_1$ we encounter at each stage of the traversal. Those events could typically say for example: there we reach a new unforeseen black vertex (forward connection) and we label its adjacent edges $c$, $c+1$, $c+2$, or there we reach a previously visited white vertex (backward connection), or there we reach an unforeseen white vertex, etc...

One can easily convince oneself that such a sequence of events, relying only on the execution march of the algorithm and not on the particular labeling of the input diagram, is in fact characteristic to the diagram. If sufficiently detailed, that sequence of events can be used to unambiguously characterize rooted trivalent diagrams. The idea now would be to consider a rooted planar tree with leaves labeled by rooted trivalent diagrams and with edges labeled by events in such a way that the sequence of events one gets along any branch from the root to a leaf is the very sequence of events that unambiguously characterize the corresponding rooted trivalent diagram.

We now obtain a usable principle of generation if we require two further properties: exhaustivity, meaning that every isomorphism class of rooted trivalent diagram gets represented on a particular leaf of the tree and non-redundancy, meaning that every such isomorphism class gets represented just once. Assuming that we spend only a constant time on each node of that tree and that the number of those nodes is linearly bounded by the number of its leaves, this would provide a constant amortized time algorithm to generate rooted trivalent diagrams.

To ease the memory requirements of the generator, we won't actually build the generation tree in memory. It will instead be realized in the calling pattern between the procedures of the generating program. Also, the program would be more useful if it generates the diagrams in permutational form instead of a sequence of events describing it. This means that we have to carry around a partial diagram that gets built while exploring the generation tree, each generating event completing that description and each backtrack reversing the particular changes we have made.

\subsection{Implementation}
\label{subsect:generatate:implementation}

The generating algorithm uses as global data two integers $c$ and $n$, a stack of integers and two integer arrays
\begin{align*}
	s_0, s_1 : \{\,1,...,n\,\} \to \{\,1,...,n\,\}
\end{align*}
representing the rooted trivalent diagram being constructed by its black and white permutations $\sigma_\bullet$ and $\sigma_\circ$ respectively. The integer $n$ represents the maximum size of the diagram being generated while the integer $c$ is the labeling counter used to attribute integer labels to the edges of rooted trivalent diagrams being generated.
The manipulation of the stack is done through the following five primitives.
\begin{align*}
	\textsc{Push}  &: \mathrm{Integer} \times \mathrm{Stack} \to \mathrm{Stack}			\\
	\textsc{Pop}  &: \mathrm{Stack} \to \mathrm{Integer} \times \mathrm{Stack}			\\
	\textsc{StackIsEmpty}  &: \mathrm{Stack} \to \mathrm{Bool}							\\
	\textsc{Mask},\textsc{Reveal}  &:\mathrm{Integer} \times \mathrm{Stack} \to \mathrm{Stack}
\end{align*}
The stack can be implemented using a doubly linked circular list represented by two zero-based arrays of integers.
\begin{align*}
	N, P : \{\,0,...,n\,\} \to \{\,0,...,n\,\}
\end{align*}
The item of index zero is just a sentinel and the stack is considered empty if the following relation holds,
\begin{align*}
N[\,0\,] = P[\,0\,] = 0.
\end{align*}
The Mask and Reveal procedures implement removal and insertion primitives using a trick popularized by Knuth \cite{knuth00} under the name of ``dancing link''. Namely, a call to the Mask procedure with parameter $s$ removes the item $s$ of the stack using the following two instructions,
\begin{align*}
	&N [\,P[\,s\,]\,] \leftarrow N[\,s\,],	\\
	&P [\,N[\,s\,]\,] \leftarrow P[\,s\,],
\end{align*}
while a subsequent call to the Reveal procedure with parameter $s$ restores the previous state of the stack, before the call to the Mask procedure, using the following two instructions,
\begin{align*}
	&N [\,P[\,s\,]\,] \leftarrow s,	\\
	&P [\,N[\,s\,]\,] \leftarrow s.
\end{align*}

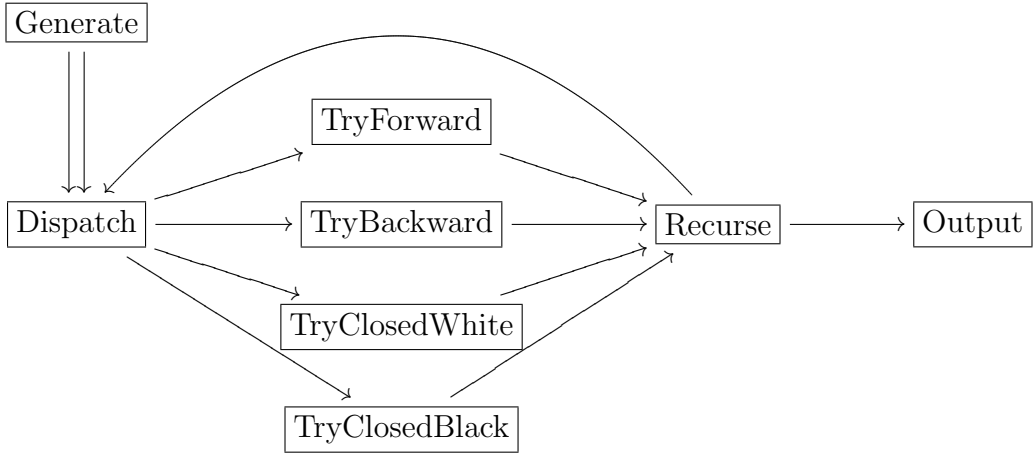
\begin{figure}
\begin{center}
\begin{align*}
\xymatrix@C=1.5cm@R=0.5cm
{
	{\fbox{\text{Generate}}}
		\ar[dd]<1mm>
		\ar[dd]<-1mm>
	\\
	&
	{\fbox{\text{TryForward}}}
		\ar[rd]
		\\
	{\fbox{\text{Dispatch}}}
		\ar[ru]
		\ar[r]
		\ar[rd]
		\ar[rdd]
	&
	{\fbox{\text{TryBackward}}}
		\ar[r]
	&
	{\fbox{\text{Recurse}}}
		\ar[r]
		\ar@/_25mm/[ll]
	&
	{\fbox{\text{Output}}}
		\\
		&
	{\fbox{\text{TryClosedWhite}}}
		\ar[ru]
		\\
	&
	{\fbox{\text{TryClosedBlack}}}
		\ar[ruu]
}
\end{align*}
\caption{Overall structure of the generating algorithm.}
\label{fig:overall:structure}
\end{center}
\end{figure}

The generating algorithm is composed of seven procedures (algorithm \ref{algo:generate} to \ref{algo:recurse}) and a user defined procedure called Output that serves as an outlet to the algorithm and that can be used for instance to do printing jobs or to collect some statistics on rooted trivalent diagrams. The overall structure of the calling pattern between those procedure is shown in figure \ref{fig:overall:structure}.
The algorithm works by a recursive exploration of the structure being constructed in a way that mimics the depth first traversal of the labeling algorithm of section 2.

The recursion has two base cases that are produced by the Generate procedure (algorithm \ref{algo:generate}) which is the main entry point of the algorithm. The inductive step of the recursion corresponds to a call to the Dispatch procedure (algorithm \ref{algo:dispatch}), which purpose is to successively handle each of the various cases one can encounter at each stage of the construction/exploration of the diagrams. This is the branching part of the generating algorithm in the sense that it is there that the generation tree forks into subtrees eventually leading to the leaves where the produced structures reside. A call to the Dispatch procedure results in a call to the Recurse procedure (algorithm \ref{algo:recurse}) through each of the four Try procedure (algorithm \ref{algo:tryforward} to \ref{algo:tryclosedblack}).
The purpose of the Recurse procedure is to call the Output procedure if the stack is empty meaning that the exploration/construction is finished and that we can thus output a finished structure or to pop an edge and call the dispatch procedure if the stack is not empty meaning that the structure being explored/constructed is not yet finished.

\subsubsection{The Generate procedure}

The Generate procedure (algorithm \ref{algo:generate}) is the main entry point of the program. It is responsible for the two base cases of the induction, namely whether the produced structure has a univalent black vertice adjacent to its root edge (handled line 3 to 5 of the procedure) or a trivalent one (handled line 7 to 13).


In the case where the vertex is univalent, the corresponding fixed point is built (line 4 of the algorithm), the labeling counter is set to 2 (the label of the next encountered edge), and then the exploration/construction continues in the direction of the white vertex by calling the Dispatch procedure with parameter 1 (line 5 of the algorithm).

In the case where the vertex is trivalent, the three edges around it are labeled 1, 2 and 3 in counterclockwise direction, the corresponding cycle in the black permutation is built by the three instructions line 8 to 10 of the algorithm and the two edges 1 and 2 are pushed on the stack for further exploration (line 11 and 12) while the edge 3 is explored in direction of its white vertex by calling the Dispatch procedure with parameter 3 (line 13 of the algorithm).
\begin{algorithm}
\dontprintsemicolon
\Begin{
	\If{$n \ge 1$}{
		$c \leftarrow 2$ \;
		$s_0[1] \leftarrow 1$ \;
		$\textsc{Dispatch }( 1 )$
	} \;
	\If{$n \ge 3$}{
		$c \leftarrow 4$ \;
		
		$s_0[1] \leftarrow 2$ \;
		$s_0[2] \leftarrow 3$ \;
		$s_0[3] \leftarrow 1$ \;
		
		$\textsc{Push }( 1 )$ \;
		$\textsc{Push }( 2 )$ \;
		$\textsc{Dispatch }( 3 )$
	} \;
}
\caption{$\textsc{Generate }(\,)$}
\label{algo:generate}
\end{algorithm}

\subsubsection{The Dispatch procedure}
\label{subsubsec:impl:dispatch}
This is the start of the induction step of the generation algorithm. The hypothesis at its start is that the two arrays $s_0$ and $s_1$ reflects the structure of a partial trivalent diagram being explored according to the depth first traversal of the labeling algorithm of section \ref{sect:char:label}.
The current edge $s$ and the labeling counter $c$ reflect the stage of the exploration. The exploration/construction is supposed to continue from the current edge $s$ in the direction of its white vertex and $c$ is the label attributed to the next unlabeled edge we encounter. At the beginning of the procedure, we don't know if that white vertex is univalent or bivalent, and if it's bivalent, we don't know which is the other edge incident to it.
There are four possible cases:

\noindent
\textbf{Case 1.} The white vertex incident to the current edge $s$ is univalent. That case is handled by a call to the TryClosedWhite procedure line 3 of the Dispatch procedure (algorithm \ref{algo:dispatch}).
We can then assume for the three other cases that this vertex is bivalent.

\noindent
\textbf{Case 2.}  The edge adjacent to the current edge by its bivalent white vertex hasn't been visited yet and the next black vertex is trivalent. This case is handled line 4 by a call to the TryForward procedure. 

\noindent
\textbf{Case 3.} As in the previous case, the adjacent edge hasn't been visited yet but here the next black vertex is univalent. This case is handled line 5 by a call to the TryClosedBlack procedure. 

\noindent
\textbf{Case 4.} The edge adjacent to the current edge by its bivalent white vertex has already been visited and thus already has a label, which we call $t$. The point is that those edges, already visited but still waiting for further exploration on their white side, are precisely those that are stored in the stack. Each edge stored in the stack corresponds to an admissible possibility for the white neighbour of the current edge. The exploration of those possibilities is done by the ``for'' statement line $6$ to $9$ of the procedure. The edge of the stack $t$ matched with the current edge $s$ is temporarily removed from the stack using the ``dancing link'' trick implemented by the Mask and Reveal primitives called line $7$ and $9$ of the procedure.

The production of new edges has to be compatible with the maximum allowed size of the diagrams $n$. That condition is checked by the two ``if'' statements line 4 and 5 of the procedure.


\begin{important}
We claim that those four cases cover all the possibilities and that they are mutually exclusive.
\end{important}

\begin{remark}
The order in which the cases are handled by the Dispatch procedure only affects the order in which the structures are produced but not the way they are labeled nor does it change the set of structure that is produced.
\end{remark}

\begin{algorithm}
\dontprintsemicolon
$\textbf{local } t : integer$  \;
\Begin{
	$\textsc{TryClosedWhite }(s)$\;
	\lIf{$c+3 \le n + 1$}{
		$\textsc{TryForward }(s)$
	} \;
	\lIf{$c+1 \le n + 1$}{
		$\textsc{TryClosedBlack }(s)$
	} \;
	
	\For{$t \in \textsc{Stack}$}{
		$\textsc{Mask }(t)$ \;
		$\textsc{TryBackward }(s,t)$ \;
		$\textsc{Reveal }(t)$ \;
	}
}
\caption{$\textsc{Dispatch }(s:integer)$}
\label{algo:dispatch}
\end{algorithm}

\subsubsection{The TryClosedWhite procedure}

It handles the case where the current edge $s$ is incident to a univalent white vertex (case 1. above), as the following picture shows.
\nopagebreak 
\begin{center}
\includegraphics{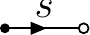}
\end{center}
It simply builds a fixed point on $s$ in the white permutation (line 2 of the procedure) then it calls the Recurse procedure.
\begin{algorithm}[h]
\dontprintsemicolon
\Begin{
	$s_1\,[\,s\,] \leftarrow s$ \;
	$\textsc{Recurse }(\,)$ \;
}
\caption{$\textsc{TryClosedWhite }(s : integer)$}
\label{algo:tryclosedwhite}
\end{algorithm}

\subsubsection{The TryForward procedure}

Its purpose is to handle the case where the current edge $s$ is incident to a bivalent white vertex followed by a trivalent black vertex (case 2. above) as the following picture shows.
\nopagebreak 
\begin{center}
\includegraphics{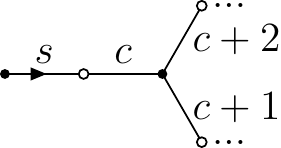}
\end{center}
The edges incident to that trivalent black vertex are supposed to never have been encountered before and are then labeled by the value $c$, $c+1$ and $c+2$ (the case where the adjacent edge has already been encountered and thus has already a label is handled by the TryBackward procedure).
Lines $2$ to $6$ build the corresponding black and white cycles in the permutation arrays.
Among the three created edges, two need further exploration on their white side so they are put on the stack by the Push instructions line $7$ and $8$. Before calling the Recurse procedure the labeling counter $c$ is increased to account for the creation of the three new edges. The sate of the stack and the labeling counter are both restored to their previous value by the instructions line $11$ to $13$, before the procedure exits. 
\begin{algorithm}[h]
\dontprintsemicolon
\Begin{
	$s_0\,[\,c\,] \leftarrow c+1$ \;
	$s_0\,[\,c+1\,] \leftarrow c+2$ \;
	$s_0\,[\,c+2\,] \leftarrow c$ \;
	$s_1\,[\,s\,] \leftarrow c$ \;
	$s_1\,[\,c\,] \leftarrow s$ \;
	$\textsc{Push }( c + 1 )$ \;
	$\textsc{Push }( c + 2 )$ \;
	$c\leftarrow c + 3$ \;
	$\textsc{Recurse } (\,)$ \;
	$c\leftarrow c -3$ \;
	$\textsc{Pop }(\,)$ \;
	$\textsc{Pop }(\,)$ \;
}
\caption{$\textsc{TryForward  }(s:integer)$}
\label{algo:tryforward}
\end{algorithm}
\subsubsection{The TryClosedBlack procedure}
It handles the case where the current edge $s$ is incident to a bivalent white vertex followed by a univalent black vertex (case 3. above), as the following picture shows.
\nopagebreak 
\begin{center}
\includegraphics{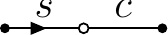}
\end{center}
It builds the white $2$-cycle and the black $1$-cycle line 2 to 4 then the labeling counter $c$ is increased line 5 to account for the creation of the new edge labeled $c$. It calls the Recurse procedures line 6 and then restores the value of the labeling counter before it exits.
\begin{algorithm}
\dontprintsemicolon
\Begin{
	$s_1\,[\,s\,] \leftarrow c$   \;
	$s_1\,[\,c\,] \leftarrow s$ \;
	$s_0\,[\,c\,] \leftarrow c$  \;
	$c\leftarrow c+1$ \;
	$\textsc{Recurse }(\,)$ \;
	$c\leftarrow c-1$ \;
}
\caption{$\textsc{TryClosedBlack }(s : integer)$}
\label{algo:tryclosedblack}
\end{algorithm}

\subsubsection{The TryBackward procedure}

It handles the case where the current edge $s$ is incident to a bivalent white vertex followed by an edge $t$ that has already been visited (case 4. above), as the following picture shows.
\nopagebreak
\begin{center}
\includegraphics{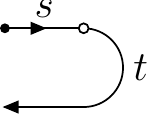}
\end{center}
The edge $t$ is chosen in the Dispatch procedure and is removed from the stack before the TryBackward procedure is called and reinstalled back after the procedure terminates. The procedure simply binds together the two edges $s$ and $t$ by building a 2-cycle in the white permutation and then calls the recurse procedure.
\begin{algorithm}[h]
\dontprintsemicolon
\Begin{
	$s_1\,[\,s\,] \leftarrow t$ \;
	$s_1\,[\,t\,] \leftarrow s$ \;
	$\textsc{Recurse }(\,)$ \;
}
\caption{$\textsc{TryBackward }(s,t : integer)$}
\label{algo:trybackward}
\end{algorithm}

\subsubsection{The Recurse procedure}

Its purpose is to check for the termination of the recursion. If the stack is empty then the recursion terminates and it calls the Output procedure. In that case, the two arrays $s_0$ and $s_1$ describe a finished rooted trivalent diagram with edges labeled from $1$, the root edge, to $c-1$, the last attributed label. The size of the diagram is thus $c-1$.
Otherwise, when the stack is not empty, an edge is poped out of the stack and the recursion continues by a call to the Dispatch procedure line 7.
\begin{algorithm}[h]
\dontprintsemicolon
$\textbf{local } k : integer$  \;
\Begin{
	\eIf{$\textsc{StackIsEmpty } (\,)$}{
		$\textsc{Output }(\,)$ \;
	}{
		$k \leftarrow \textsc{Pop } (\,)$ \;
		$\textsc{Dispatch } ( k )$ \;
		$\textsc{Push } ( k )$ \;
	}
}
\caption{$\textsc{Recurse }(\,)$}
\label{algo:recurse}
\end{algorithm}

\subsection{Correctness}

In this section we show that the algorithm is correct and provide a complexity analysis of its execution time. This complexity analysis is based on a study of the structure of the execution tree of the algorithm and rely on the assumption that it is finite. We first prove that assumption.

\begin{lemma}
The generating algorithm terminates in a finite amount of time.
\end{lemma}

\begin{demonstration}
Looking at the procedures we see that each of them takes only a finite time to complete provided the procedure that are called also take in turn a finite time to complete. So the proof reduces to show that only a finite number of procedures is called, meaning that the execution tree of the algorithm is finite.
Using Knig lemma on trees one has to show that all the branches of that tree are finite. We show that there is a non-negative integer quantity that is strictly decreasing along every branches.
If we denote by $n_s$ the size of the stack then $\mu = 2\,n_s + n - c  + 1$ is such a quantity, where $n$ is the maximum size of the diagrams being generated and $c$ is the labeling counter of the algorithm.
It is non-negative because $n_s$ is non-negative and because according to the ``if'' conditions of the Dispatch procedure, $n - c + 1$ is also non-negative, meaning that no label exceeding $n$ is ever attributed to an edge. In the following table we summarize the changes in the value of $n_s$, $c$ and $\mu$ after a cycle Dispatch $\rightarrow$ Try-- $\rightarrow$ Recurse $\rightarrow$ Dispatch has been completed in each of the four cases described in section \ref{subsubsec:impl:dispatch} (\confer figure \ref{fig:overall:structure}).
\begin{center}
\begin{tabular}{|ccc|ccc|ccc|ccc|ccc|}
\hline
&&&& case 1. &&& case 2. &&& case 3. &&& case 4. & \\
\hline
&$n'_s$ &&& $n_s - 1$ &&& $n_s +1$ &&& $n_s-1$ &&& $n_s-2$ &  \\
&$c'$  &&& $c$ &&& $c+3$ &&& $c+1$ &&& $c$  & \\
&$\mu'$ &&& $\mu-2$ &&& $\mu-1$ &&& $\mu-3$ &&& $\mu-4$ &  \\
\hline
\end{tabular}
\end{center}
In each case, $\mu' < \mu$ so the quantity $\mu$ is strictly decreasing along every branch of the execution tree.
\hfill $\Box$
\end{demonstration}

\begin{lemma}
The rooted trivalent diagrams produced by the generating algorithm are labeled according to the characteristic labeling of section \ref{sect:char:label}.
\end{lemma}

\begin{demonstration}
The proof is by induction.
Assuming that the stack, the labeling counter and the labels of the already generated cycles agree with the corresponding state of the labeling procedure of section \ref{sect:char:label} at the start of a call to the Dispatch procedure (induction hypothesis) one can check that the execution of the algrithm through each of the four Try-- procedures (algorithm \ref{algo:tryforward} to \ref{algo:tryclosedblack}) each preserves that hypothesis, that is the new cycles introduced in the permutations $s_0$ and $s_1$ are labeled consistently with that of the labeling procedure and that the state of the stack also match the one found in the labeling after visiting those new cycles.
Finally, one has to check that the two base cases produced by the Generate procedure are also consistent with the characteristic labeling. This is immediate and completes the induction.
\hfill $\Box$
\end{demonstration}

\begin{theorem}
The generating algorithm produces an exhaustive and non-redundant list of rooted trivalent diagrams.
\end{theorem}

\begin{demonstration}
Exhaustivity comes primarily by induction from the local exhaustivity claim of the case analysis of section \ref{subsubsec:impl:dispatch}. The mutual exclusion of the cases ensures that different rooted diagrams are produced, at least differing in the way they are labeled (the labeling counter is strictly increasing during the generation process of a structure), but since the structures are produced in their characteristic labeling, none of them could be isomorphic.
\hfill $\Box$
\end{demonstration}

\subsection{Average Time Complexity}

In this section we prove the main property of the algorithm, that it spends a constant amortized time generating each structure. One way to do the proof could be to express an estimate of the total execution time and the number of structure produced and show that the quotient of those two quantity is bounded independently of the size of the produced structure.
We propose instead a proof of the majoration based on the following principle and a careful analysis of the execution tree of the generating algorithm.

\theoremstyle{plain}
\newtheorem*{balanceprinciple}{Balance Principle}
\begin{balanceprinciple}
In a finite tree, the number of leaves is greater than the number of its node having degree at least $2$.
\end{balanceprinciple}

\begin{demonstration}
The proof is by an easy recurrence on the number of internal node. Every finite trees can be constructed starting from a one node tree by successive replacement of a leaf by an internal node having only leaves as sons.
The starting tree satisfies the balance principle so the recurrence is initialized. Assuming by recurrence that a finite tree satisfies the principle, and replacing one of its leaves by a internal node having $k \ge 1$ leaves as sons, one increases the number of internal nodes by 1 and the number of leaves by $k-1 \ge 0$. If $k \ge 2$ the number of nodes having degree at least 2 is increased by 1, but the number of leaves is also increased by $k-1\ge 1$ so the resulting tree still satisfies the principle. The recurrence is complete.
\hfill $\Box$
\end{demonstration}

\begin{lemma}
The total execution time of the generating algorithm is $O(C_n)$, where $C_n$ is the number of procedure called during the execution.
\end{lemma}

\begin{demonstration}
The only procedure of the algorithm that contains a loop an thus can have an arbitrary long execution time is the Dispatch procedure. Since each iteration of the loop has a constant execution time and that each time the TryBackward procedure is called, we can transfer the expanse of the iteration to the TryBackward procedure and then assume the Dispatch procedure to have a bounded execution time.
This way every procedure is assumed to have a bounded execution time so that the total execution time is proportional to $C_n$.
\hfill $\Box$
\end{demonstration}

Let $D_n$, $R_n$, $O_n$ denote respectively the total number of time the Dispatch, Recurse and Output procedures are called and let $T_n$ denote the total number of time one of the four Try-- procedures is called, $n$ is the maximum size of the structures being produced.

\begin{lemma}
We have $D_n \le 2\,O_n$.
\end{lemma}

\begin{demonstration}
Let $D'_n$ denote the number of call to the Dispatch procedure that have an out degree at most $2$. The leaves of the execution tree of the algorithm are calls to the Output procedure because the other procedures all have out degree at least  1, hence by the balance principle above $D'_n \le O_n$. The Dispatch procedure has an out degree at least one because it calls the TryClosedWhite procedure unconditionally line $3$, and when its outdegree is $1$ then it means that the stack is empty (line $6$) and that $c \le n$ (line 4 and 5). This means that all the $n$ edges has been labeled and that the call to the Recurse procedure subsequent to the call to the TryClosedWhite will result in a call to the Output procedure and no further call to the Dispatch procedure. Therefore we see that the Dispatch procedure can have an out degree of 1 but that can only happend once at the end of each branch. We thus have $D_n \le D'_n + O_n$ and then $D_n \le 2\,O_n$.
\hfill $\Box$
\end{demonstration}

\begin{theorem}[CAT property]
The average time spent by the generating algorithm to produce each structure is bounded independently of their size.
\end{theorem}
\begin{demonstration}
Let $C_n$ denote the total number of procedure called during the execution of the algorithm. The number of structure produced by the algorithm is equal to $O_n$. Since the total execution time of the algorithm is $O(C_n)$ the average time spent producing each structure is $O(C_n/O_n)$. We have to show that the quotient $C_n/O_n$ is bounded independently of $n$. We clearly have $C_n = 1+D_n+R_n+T_n+O_n$, the $+1$ accounts for the first call to the Generate procedure. Since the Recurse procedure is only called by one of the four Try-- procedures and that each one calls it exactly once we have $T_n = R_n$. Since the Recurse procedure is called twice by the Generate procedure and that the Recurse procedure calls one of the Recurse or Output procedures, we have $R_n + 2 = D_n +  O_n$.  Therefore we have,
\begin{align*}
	&& C_n &=  1+D_n+R_n+T_n+O_n\\
	&&	&= 1+D_n+2\,R_n + O_n 		&& \text{as $T_n = R_n$,} \\
	&&	&= 3\,D_n + 3\,O_n - 3 		&& \text{as $R_n + 2 = D_n +  O_n$,}\\
	&&	&\le 9\, O_n				&& \text{as $D_n \le 2\,O_n$,}
\end{align*}
and then $C_n/O_n \le 9$. The bound on the quotient is independent of $n$ as announced.
\hfill $\Box$
\end{demonstration}

\section[First Application: Modular Group and Unrooted Trivalent Diagrams]{First Application: Modular Group and\\ Unrooted Trivalent Diagrams}
\label{sect:application:modular:group}

We recall that the modular group $\mathrm{PSL}_2(\ZZ)$ is the following group of $2$ by $2$ integer matrices with unit determinant:
\begin{align*}
	\mathrm{PSL}_2(\ZZ) = \left\{\,\pm\begin{pmatrix} a & b \\ c & d \end{pmatrix} \in \mathcal{M}_2(\ZZ)/\pm\mathrm{Id}\,\bigg|\,ad-bc = 1\,\right\}.
\end{align*}
There are many possible finite presentations for this group and we shall stick to the following:
\begin{align*}
	\mathrm{PSL}_2(\ZZ) = \langle\, A, B \,|\, A^2 = B^3 = 1 \,\rangle
\end{align*}
with $A$ and $B$ being the following two matrices:
\begin{align*}
	A = \pm\left(\begin{array}{rr} 0 & -1 \\ 1 & 0 \end{array}\right)
	\quad\text{ and } \quad
	B = \pm\left(\begin{array}{rr} 0  & 1 \\ -1 & 1 \end{array}\right),
\end{align*}
for it renders explicit the following isomorphism:
\begin{align*}
	\mathrm{PSL}_2(\ZZ) \simeq \ZZ/2\ZZ \ast \ZZ/3\ZZ.
\end{align*}



\subsection{Displacements Groups}

The modular group acts naturally on the set of edges of any trivalent diagram. This action is generated by the two \emph{elementary moves},
\begin{align*}
	a \cdot A = \sigma_\circ (a)
	\quad\text{ and }\quad
	a \cdot B = \sigma_\bullet(a).
\end{align*}
The elementary move $A$ acts by exchanging positions of the two adjacent edges of any bivalent white vertex and by fixing the only adjacent edge of any univalent white vertex. Similarly, the elementary move $B$ acts by cyclically permuting the three edges incident to any trivalent black vertex and by fixing the only edge incident to any univalent black vertex.
\begin{figure}[h]
\begin{center}
\includegraphics{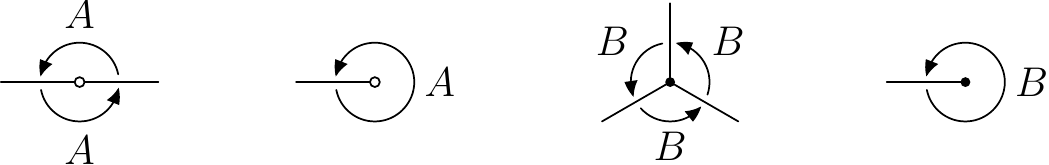}
\caption{Here is a picture of the result of the action of the two elementary moves $A$ and $B$ on the various sorts of edges.}
\label{fig:modular:action}
\end{center}
\end{figure}

Given any trivalent diagram $\Gr$, the two elementary moves just described generate a group $\Phi_\Gr$ called the displacement group of $\Gr$. It is easily verified that it is the quotient group of $\mathrm{PSL}_2(\ZZ)$ by the kernel of the group action $\rho : \mathrm{PSL}_2(\ZZ) \to \mathfrak{S}_{\Gr_{-}}$.
The modular group has therefore a universal status with respect to that construction, it can be considered as the universal group of displacements for the species of trivalent diagrams. If one restricts one's attention to finite trivalent diagrams, the profinite completion of $\mathrm{PSL}_2(\ZZ)$ would be a more appropriate candidate for that purpose.

\subsection{Unrooted Planar Binary Trees}
\label{sect:unrooted:planar:bin:tree}

One can associate to any (unrooted) planar binary tree $\Theta$ a connected and acyclic trivalent diagram $\Gr$, called its \emph{enriched barycentric subdivision} $\Gr = \Theta^{sb+}$, by putting an extra white vertex in the middle of every edges of $\Theta$. The set of directed edge of $\Theta$ and that of undirected edges of $\Gr$ are in bijection in two natural ways.

\begin{figure}
\begin{center}
\includegraphics{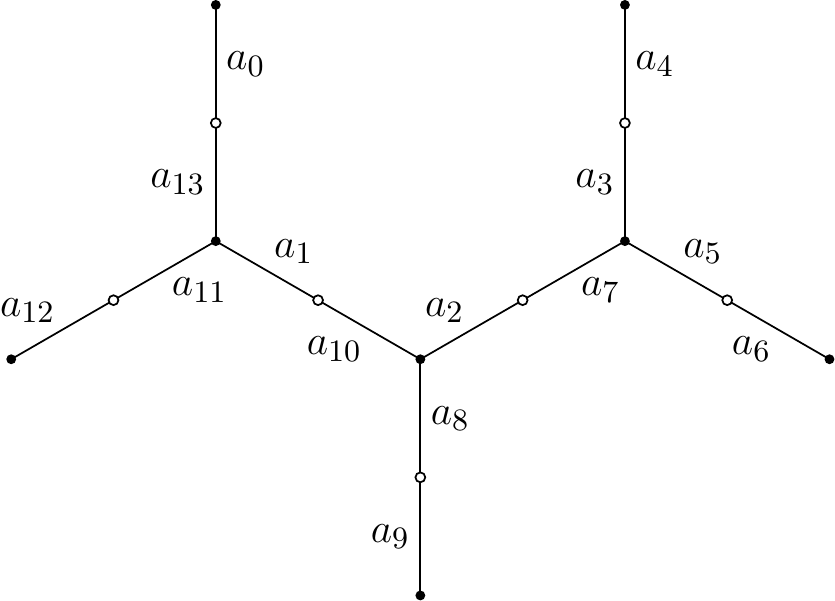}
\caption{We see in this example the result of iterating the elementary move $T$ on the edges of a binary tree. The edges are labeled by $a_k$ where $a_{k+1} = a_{k} \cdot T$. This can be used to implement depth-first traversals in a purely iterative way.}
\label{fig:depth:first:bin:tree}
\end{center}
\end{figure}

There is another famous presentation of the modular group. It is given by two generators $S$ and $T$ and two relations as follows,
\begin{align*}
	\mathrm{PSL}_2(\ZZ) = \langle\, S, T \,|\, S^2 = (ST)^3 = 1 \,\rangle .
\end{align*}
with $S$ and $T$ being the following two matrices:
\begin{align*}
	S = \pm\left(\begin{array}{rr} 0 & -1 \\ 1 & 0 \end{array}\right)
	\quad\text{ and } \quad
	T = \pm\left(\begin{array}{rr} 1  & 1 \\ 0 & 1\end{array}\right).
\end{align*}
The conversion between the two presentations is done through the application of the following rules:
\begin{align*}
	A &\to S, &
	S & \to A, \\
	B &\to (ST)^2, &
	T &\to AB\inv.
\end{align*}

Here are two basic criteria relating connectedness and acyclicity of finite trivalent diagrams to the transitivity of the action of some displacement groups:
\begin{enumerate}[\quad 1)\:\:]
\item A finite trivalent diagram $\Gr$ is \emph{connected} if and only if its displacement group $\Phi_\Gr$ acts transitively on its set of edges.
\\
\item If it is a tree, then the subgroup $\Psi_\Gr$ of its displacements generated by the elementary move $T$ acts transitively on its set of edges.
\end{enumerate}

There is a natural bijection between trivalent diagrams having no univalent white vertex and those having no univalent black vertex. It works by removing every univalent black vertex and the adjacent edges in one direction and by growing every univalent white vertex with a new edge and a new univalent black vertex in the other direction.
This bijection is compatible with connectedness and acyclicity and thus restricts from the class of trivalent diagrams to that of planar binary trees, and we recover the classical bijection between complete and incomplete planar binary trees.

\subsection{Classification Principle}

To any connected rooted trivalent diagram, one can moreover associate the subgroup of the $\mathrm{PSL}_2(\ZZ)$ consisting of elements that fix the distinguished edge of the diagram. We proved in \cite{vidal06} that this correspondence is one to one and we presented an enumeration in the form of a generating series that agrees perfectly with the number of structures generated by algorithm \ref{algo:generate}.

\begin{table}
\begin{center}
\includegraphics{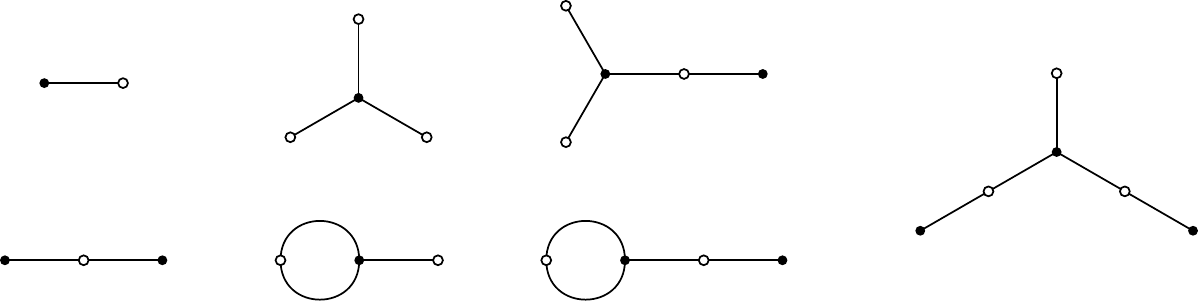}
\caption{Trivalent diagrams of size up to five, up to isomorphism.}
\label{tab:diag:3:taille:inf:5}
\end{center}
\end{table}

\begin{table}
\begin{center}
\includegraphics{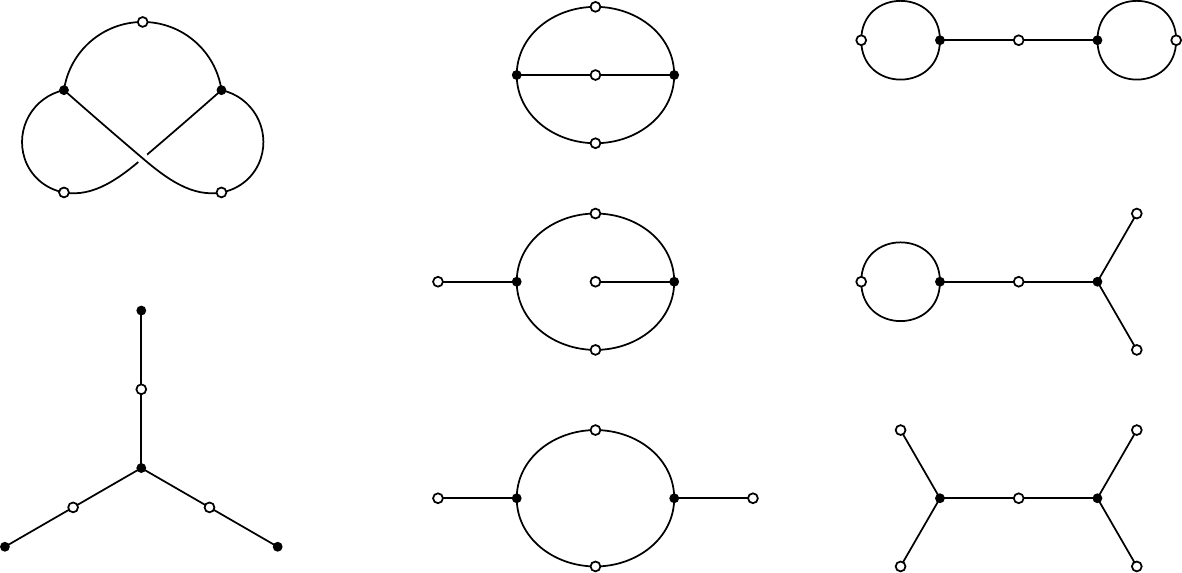}
\caption{Trivalent diagrams of size six, up to isomorphism.}
\label{tab:diag:3:taille:6}
\end{center}
\end{table}

\begin{table}
\begin{center}
\includegraphics{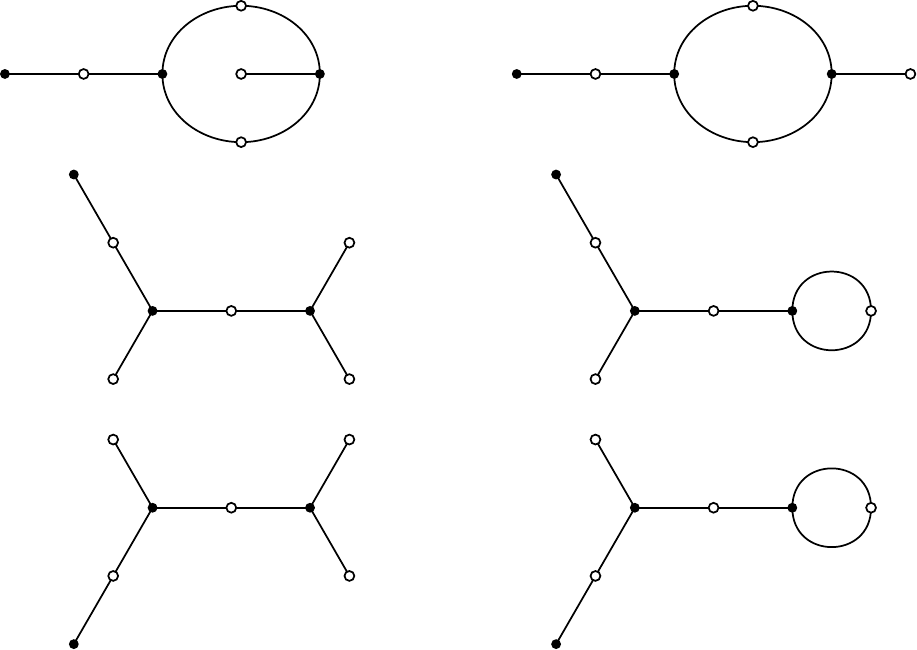}
\caption{Trivalent diagrams of size seven, up to isomorphism.}
\label{tab:diag:3:taille:7}
\end{center}
\end{table}

\begin{table}
\begin{center}
\includegraphics{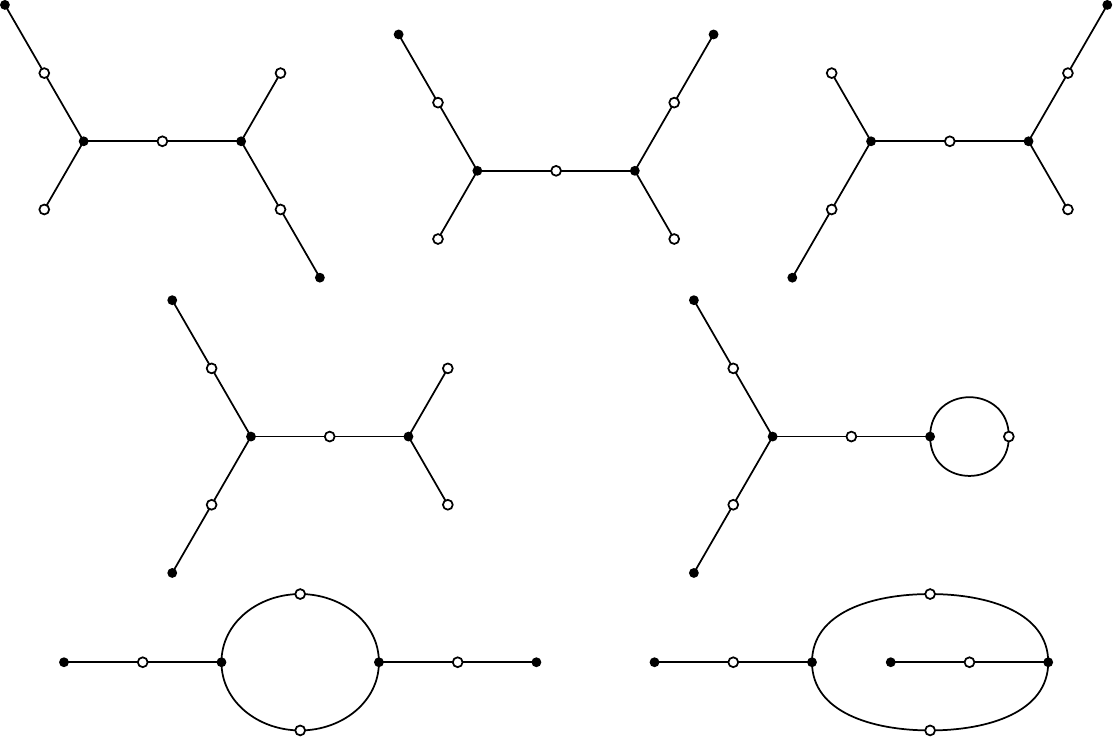}
\caption{Trivalent diagrams of size eight, up to isomorphism.}
\label{tab:diag:3:taille:8}
\end{center}
\end{table}

\begin{table}
\begin{center}
\includegraphics{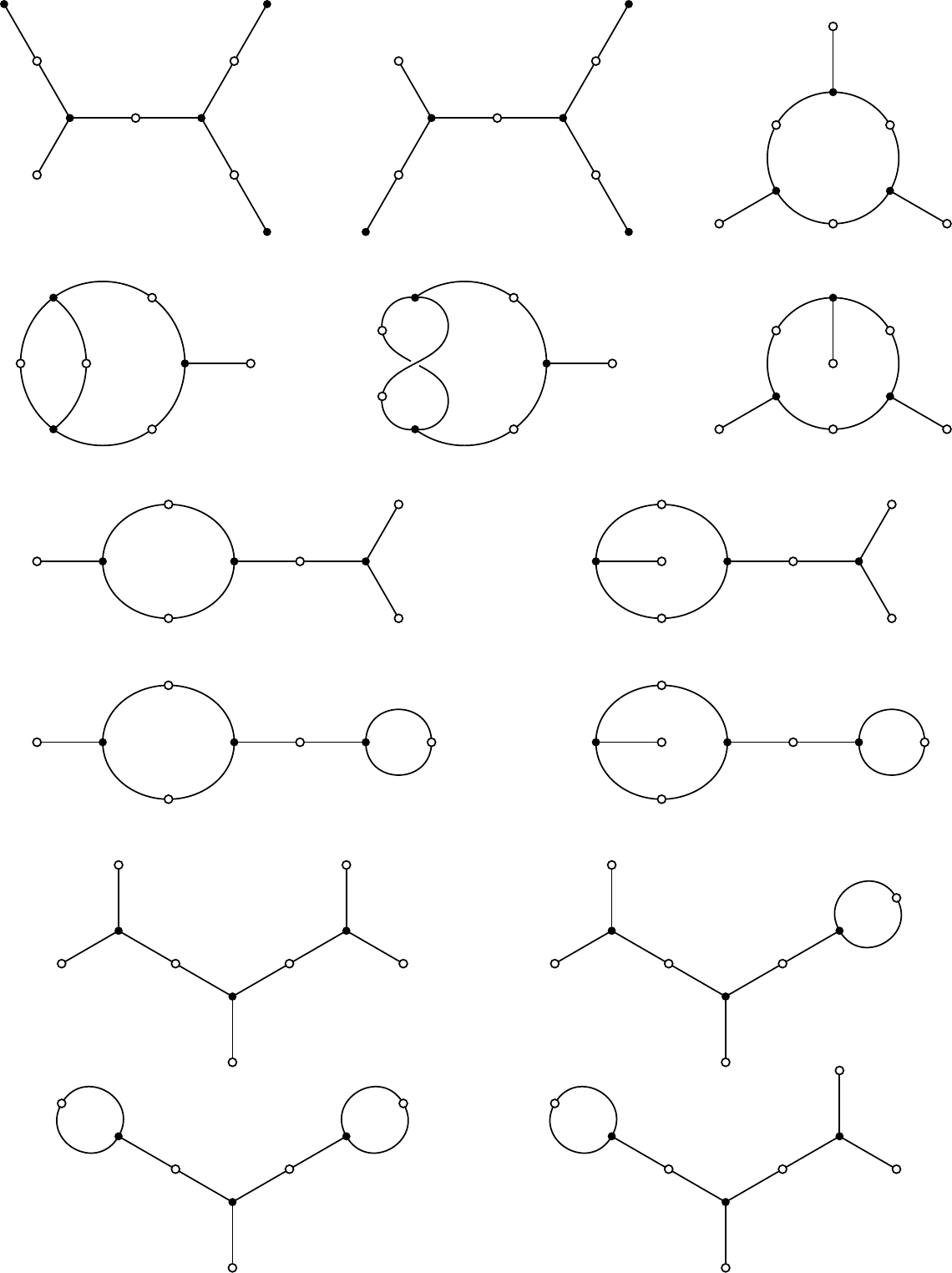}
\caption{Trivalent diagrams of size nine, up to isomorphism.}
\label{tab:diag:3:taille:9}
\end{center}
\end{table}

Now, if one changes the distinguished edge of a rooted trivalent diagram, the corresponding subgroups get conjugated. We have moreover proved that two subgroups in the modular group are conjugated if and only if the associated rooted trivalent diagrams only differ by the position of their distinguished edges. It follows that the unrooted trivalent diagrams correspond in a one-to-one fashion to the conjuguacy classes of subgroups of the modular group.

We gave in \cite{vidal06} the exhaustive list of trivalent diagrams of size up to nine. That list was computed by hand in a non-systematic  fashion. One intent of Algorithm \ref{algo:generate} of section \ref{subsect:generatate:implementation} is to permit  a retroactive validation of both those tables and the associated generating series that we reproduce here in tables \ref{tab:ser:gen:d3dot:unlabeled} and \ref{tab:serie:diag:triv}~; they are parts of the online encyclopedia of integer sequences \cite{sloan:oeis} under the references (A005133) and (A121350).

We generated unrooted trivalent diagrams using a method already described in \cite{walsh83}: by imposing a linear order on the set of rooted trivalent diagrams with a given number of edges, generating them all and accepting only those that are maximal, according to the linear order, among all those that differ only in the position of their root. Table \ref{tab:diag:3:taille:inf:5} to \ref{tab:diag:3:taille:9} show all the unrooted trivalent diagrams with up to 9 edges.

\begin{table}
\begin{center}
\begin{align*}
\tilde{D}_3^\bullet(t)\, &= \,t+{t}^{2}+4\,{t}^{3}+8\,{t}^{4}+5\,{t}^{5}+22\,{t}^{6}+42\,{t}^{7}+40
\,{t}^{8}+120\,{t}^{9}+265\,{t}^{10}+286\,{t}^{11} \\
&\quad+764\,{t}^{12}+1729
\,{t}^{13}+2198\,{t}^{14}+5168\,{t}^{15}+12144\,{t}^{16}+17034\,{t}^{
17}+37702\,{t}^{18} \\
&\quad+88958\,{t}^{19}+136584\,{t}^{20}+288270\,{t}^{21}+
682572\,{t}^{22}+1118996\,{t}^{23} \\
&\quad+2306464\,{t}^{24}+5428800\,{t}^{25}
+9409517\,{t}^{26}+19103988\,{t}^{27}+44701696\,{t}^{28} \\
&\quad+80904113\,{t}
^{29}+163344502\,{t}^{30}+379249288\,{t}^{31}+711598944\,{t}^{32} \\
&\quad+
1434840718\,{t}^{33}+3308997062\,{t}^{34}+6391673638\,{t}^{35}+
12921383032\,{t}^{36} \\
&\quad+29611074174\,{t}^{37}+58602591708\,{t}^{38}+
119001063028\,{t}^{39} \\
&\quad+271331133136\,{t}^{40}+547872065136\,{t}^{41}+
1119204224666\,{t}^{42} \\
&\quad+2541384297716\,{t}^{43}+5219606253184\,{t}^{44
}+10733985041978\,{t}^{45} \\
&\quad+24300914061436\,{t}^{46}+50635071045768\,{t
}^{47}+104875736986272\,{t}^{48} \\
&\quad+236934212877684\,{t}^{49}+
499877970985660\,{t}^{50}+o(t^{50})
\end{align*}
\caption {Order fifty development of the generating series $\tilde{D}^\bullet_3(t)$ giving as the coefficient of $t^n$ the number of connected rooted trivalent diagrams with $n$ edges (A005133) which is also the number of index $n$ subgroups in the modular group $\mathrm{PSL}_2(\ZZ)$.}
\label{tab:ser:gen:d3dot:unlabeled}
\end{center}
\end{table}

\begin{table}
\begin{center}
\begin{equation*}
\begin{split}
\tilde{D}_3(t)\, &= \,t+{t}^{2}+2\,{t}^{3}+2\,{t}^{4}+{t}^{5}+8\,{t}^{6}+6\,{t}^{7}+7\,{t}^{8}+14\,{t}^{9}+27\,{t}^{10}+26\,{t}^{11}\\
&\quad+80\,{t}^{12}
+133\,{t}^{13}+170\,{t}^{14}+348\,{t}^{15}+765\,{t}^{16}+1002\,{t}^{17}+2176\,{t}^{18}\\
&\quad+4682\,{t}^{19}
+6931\,{t}^{20}+13740\,{t}^{21}+31085\,{t}^{22}+48652\,{t}^{23}+96682\,{t}^{24}\\
&\quad+217152\,{t}^{25}
+362779\,{t}^{26}+707590\,{t}^{27}+1597130\,{t}^{28}+2789797\,{t}^{29}\\
&\quad+5449439\,{t}^{30}
+12233848\,{t}^{31}+22245655\,{t}^{32}+43480188\,{t}^{33}\\
&\quad+97330468\,{t}^{34}+182619250\,{t}^{35}+358968639\,{t}^{36}+800299302\,{t}^{37} \\
&\quad+1542254973\,{t}^{38}+3051310056\,{t}^{39}+6783358130\,{t}^{40}+13362733296\,{t}^{41}\\
&\quad+26648120027\,{t}^{42} +59101960412\,{t}^{43}+118628268978\,{t}^{44}\\
&\quad+238533003938\,{t}^{45}+528281671324\,{t}^{46}+1077341937144\,{t}^{47}\\
&\quad+2184915316390\,{t}^{48}+4835392099548\,{t}^{49}+9997568771074\,{t}^{50}
 + o(t^{50})
\end{split}
\end{equation*}
\caption{Order fifty development of the generating series $\tilde{D}_3(t)$ giving as the coefficient of $t^n$ the number of connected unrooted trivalent diagrams with $n$ edges (A121350) which is also the number of conjugacy classes of index $n$ subgroup in the modular group $\mathrm{PSL}_2(\ZZ)$.}
\label{tab:serie:diag:triv}
\end{center}
\end{table}

\section{Second Application: Triangular Maps}
\label{sect:application:triangular:maps}


By an (\emph{oriented}) \emph{triangular map} we mean a finite polyhedral structure composed of vertex, directed edges, and oriented triangular faces with an incidence relation between them.
Suggestively enough, a directed edge, also called an \emph{arc}, is bordered by an ordered pair of vertices, which we call its \emph{origin} and its \emph{destination} respectively,
such that triangular faces are each bordered by a cycle of three arcs whose destination coincides with the origin of the following arc in cyclic order.

The following definition is useful in grasping the incidence relations of a combinatorial map but insufficient because it lacks some traversal information such as the cyclic orientation of the faces. 

\begin{definition}
A \emph{combinatorial pre-map} $\Map$ is given by three sets $\Map_0$ $\Map_1$ and $\Map_2$ and five mappings $s,t : \Map_1 \to \Map_0$, $\ell, r : \Map_1 \to \Map_2$ and $.\inv : \Map_1 \to \Map_1$
satisfying the following conditions for all elements $a$ of $\Map_1$,
\begin{align*}
	s(a\inv) &= t(a)		&\ell(a\inv) &= r(a)		&(a\inv)\inv &= a \\
	t(a\inv) &= s(a)		&r(a\inv) &= \ell(a)		&a\inv &\neq a
\end{align*}
the four mappings $s$, $t$, $\ell$ and $r$ are further assumed to be \emph{surjective}.
\end{definition}

The elements of the three sets $\Map_0$, $\Map_1$, and $\Map_2$ are the \emph{vertices}, the \emph{arcs} (directed edges), and \emph{faces} of the combinatorial map respectively. The two mappings $s$ and $t$ map any arc $a$ to its \emph{origin} $s(a)$ and \emph{destination} $t(a)$. The two mappings $\ell$ and $r$ map any arc $a$ to its \emph{left-hand face} $\ell(a)$ and \emph{right-hand face} $r(a)$. Finally, the mapping $.\inv$ maps any arc $a$ to its \emph{inverse} $a\inv$ obtained by reversing its direction.

\begin{definition}
A \emph{morphism}
$\phi$ between two combinatorial pre-maps $\Map$ and $\Map'$ is a triple of mappings
$\phi_0 : \Map_0 \to \Map_0'$, $\phi_1 : \Map_1 \to \Map_1'$ and $\phi_2 : \Map_2 \to \Map_2'$ compatible with the five structure mappings in the sense that the following diagrams are commutatives.
\begin{align*}
\xymatrix@C=1.5cm@R=1.5cm
{
	{\Map_1}
		\ar[r]^{\phi_1}
		\ar[d]_{s,t}			&
	{\Map_1'}
		\ar[d]^{s,t}			\\
	{\Map_0}
		\ar[r]_{\phi_0}		&
	{\Map_0'}
}
&&
\xymatrix@C=1.5cm@R=1.5cm
{
	{\Map_1}
		\ar[r]^{\phi_1}
		\ar[d]_{\ell,r}			&
	{\Map_1'}
		\ar[d]^{\ell,r}			\\
	{\Map_2}
		\ar[r]_{\phi_2}		&
	{\Map_2'}
}
&&
\xymatrix@C=1.5cm@R=1.5cm
{
	{\Map_1}
		\ar[r]^{\phi_1}
		\ar[d]_{.\inv}	&
	{\Map_1'}
		\ar[d]^{.\inv}	\\
	{\Map_1}
		\ar[r]_{\phi_1}		&
	{\Map_1'}
}
\end{align*}
When the three mappings are bijections the morphisme is an \emph{isomorphism}.
\end{definition}

\subsection{Cyclic Orientation}

In a given combinatorial pre-map $\Map$, the \emph{inner border} of a face $f$ is the set $\ell\inv(f) = \{\,a\in \Map_1\,|\,\ell(a)=f\,\}$ of arcs having $f$ as their left-hand face.
A combinatorial pre-map is said to be \emph{strictly triangular} if each of its faces has exactly three arcs in its inner border.
The following definition describes the traversal information lacking to a triangular combinatorial pre-map to fully describe a triangular map.

\begin{definition}
A (\emph{strictly}) \emph{triangular map} $\Map$ is a strictly triangular combinatorial pre-map together with a cyclic order imposed on the three arcs that have the same left-hand face: $a+1$ denotes the next arc after $a$ in this cyclic order.
\end{definition}

\begin{definition}
A \emph{morphism} $\phi$ (respectively, an \emph{isomorphism}) between two triangular maps $\Map$ and $\Map'$ is a morphism (respectively, an isomorphism) of the underlying combinatorial pre-maps 
that preserves the cyclic order of the previous definition.
\end{definition}

\subsection{Associated Trivalent Diagram}

To any triangular map $\Map$ one can associate in a natural way a trivalent diagram $\Gr = \Map^{inc}$ which is called the \emph{incidence diagram} of the triangular map.
To any face of $\Map$ one associate a black vertex of $\Gr$ to any edge of $\Map$ one associate a white vertex of $\Gr$ then one connect black and white vertices of $\Gr$ according to the incidence of their $\Map$ counterpart. The cyclic orientations around black vertice is the one given by the cyclic orientation of the corresponding triangular faces of $\Map$. Table \ref{tab:map:2} gives three simple examples of the correspondance. This operation has already been described by Walsh in \cite{walsh75}, namely: taking the face-vertex dual of $\Map$ one gets a trivalent map and then applying the construction of Walsh, one gets the incidence diagram $\Gr$.

The \emph{incidence diagram} of a triangular map $\Map$ is the trivalent diagram, denoted $\Map^{inc}$, which sets of  white vertices, black vertices and edges are the following,
\begin{align*}
	\Map^{inc}_\circ &= \{\fra_a\}_{a\in\Map_1^*}				&
	\Map^{inc}_\bullet &= \{\frb_f\}_{f\in\Map_2}				&
	\Map^{inc}_{-} &= \{\frc_a\}_{a\in\Map_1}
\end{align*}
where $\Map_1^*$ is the set of undirected edges of $\Map$,
and which structure mappings $\partial_\circ : \Map^{inc}_{-} \to \Map^{inc}_\circ$, $\partial_\bullet : \Map^{inc}_{-} \to \Map^{inc}_\bullet$ and the permutations $\sigma_\bullet$ and $\sigma_\circ$ of $ \Map^{inc}_{-}$ are defined by the following relations,
\begin{align*}
	\partial_\circ ( \frc_a ) &= \fra_{\pi(a)}	,					&
	\sigma_\circ(\frc_a) &= \frc_{a\inv},						\\
	\partial_\bullet ( \frc_a ) &= \frb_{\ell(a)} ,					&
	\sigma_\bullet(\frc_a) &= \frc_{a+1},
\end{align*}
where $\pi(a)$ is the the undirected version of the arc $a$.
This operation is functorial for it is easily extended to morphisms of triangular maps by the following process: to any morphism $\phi$ between two triangular maps $\Map$ and $\Map'$, we associate a morphism denoted $\phi^{inc}$ between the corresponding incidence diagrams $\Map^{inc}$ and $(\Map')^{inc}$. It is given by the following relations,
\begin{align*}
	\phi^{inc}_\circ(\fra_a) &= \fra'_{\phi_1^*(a)}				&
	\phi^{inc}_\bullet(\frb_f) &= \frb'_{\phi_2(f)}				&
	\phi^{inc}_{-}(\frc_a) &= \frc'_{\phi_1(a)}
\end{align*}
Functoriality should be obvious by careful inspection.

The $.^{inc}$ functor we get by what precedes is full and faithful, but not essentially surjective because trivalent diagrams we get as the incidence diagram of a trivalent map $\Map$ have no univalent white vertex nor univalent black vertex. Let's call \emph{regular} a trivalent diagram having no univalent vertex, the $.^{inc}$ functor is essentially surjective on the full subcategory of regular trivalent diagrams, and so,

\begin{theorem}
The $.^{inc}$ functor realize an equivalence between the category of triangular maps and the full subcategory of regular trivalent diagrams.
\end{theorem}

\begin{demonstration}
To prove this theorem we shall describe a \emph{reconstruction} operation, which associate to any regular trivalent diagram $\Gr$ a triangular map denoted $\Gr^{map}$, with functorial property and such that for all regular trivalent diagrams $\Gr$ and triangular maps $\Map$, one have two natural \emph{reciprocity isomorphisms} as follows,
\begin{align*}
	(\Gr^{map})^{inc} \underset{nat.}{\simeq} \Gr
	\quad\text{ and }\quad
	(\Map^{inc})^{map} \underset{nat.}{\simeq} \Map
\end{align*}
We shall first introduce some notations. Let's call $\Psi_\Gr$ the subgroup of $\Phi_\Gr$ generated by the elementary move $T$ (\confer section \ref{sect:unrooted:planar:bin:tree}) and lets call $\pi : \Gr_{-} \to \Gr_{-}/\Psi_\Gr$ the natural projection. The mapping induced between $\Gr_{-}/\Psi_\Gr$ and $\Gr'_{-}/\Psi_\Gr$ by an equivariant mapping $\phi : \Gr_{-} \to \Gr'_{-}$ will be denoted $\phi_{\Psi}$.
The sets of vertices, edges and faces of the reconstructed map are the following,
\begin{align*}
	\Gr^{map}_0 &= \{\,\frd_x\,\}_{x\in\Gr_{-}/\Psi_\Gr}			&
	\Gr^{map}_1 &= \{\,\fre_a\,\}_{a\in\Gr_{-}}					&
	\Gr^{map}_2 &= \{\,\frf_y\,\}_{y\in\Gr_\bullet}
\end{align*}
The five structure mappings $s,t : \Gr^{map}_1 \to \Gr^{map}_0$,  $r,\ell : \Gr^{map}_1 \to \Gr^{map}_2$ and $.\inv : \Gr^{map}_1 \to \Gr^{map}_1$ and the group action $+:\Gr^{map}_1 \times \ZZ \to \Gr^{map}_1$ of the reconstructed map are given by the following equations,
\begin{align*}
	s (\fre_a) &= \frd_{\pi(a)}								&
	\ell (\fre_a) &= \frf_{\partial_\bullet(a)}					&
	\fre_a\inv &= \fre_{\sigma_\circ(a)}							\\
	t(\fre_a) &= \frd_{\pi(\sigma_\circ(a))}						&
	r(\fre_a) &= \frf_{\partial_\bullet(a\inv)}					&
	\fre_a+1 &= \fre_{\sigma_\bullet(a)}
\end{align*}
The construction then extends to morphisms in the sense that any morphism $\phi$ between to regular trivalent diagrams $\Gr$ and $\Gr'$ induces a morphism $\phi^{map}$ between the two reconstructed maps $\Gr^{map}$ and $(\Gr')^{map}$ which three components are the following,
\begin{align*}
	\phi^{map}_0(\frd_x) &= \frd_{\phi_{-,\Psi}(x)}					&
	\phi^{map}_1(\fre_a) &= \fre_{\phi_{-}(a)}						&
	\phi^{map}_2(\frf_y) &= \frf_{\phi_\bullet(y)}
\end{align*}
The fonctoriality of the reconstruction operation and the two reciprocity isomorphisms should be clear by careful inspection.
\hfill$\Box$
\end{demonstration}

\begin{remark}
The content of the above theorem is nothing but a specific notion of \emph{Poincar duality}.
\end{remark}

\subsection{Exhaustive Generation of Triangular Maps}

\begin{table}
\begin{center}
\begin{tabular}{ccc}
\begin{tabular}{c}
\includegraphics{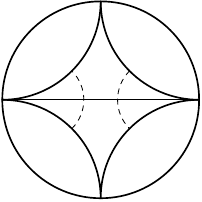} 
\end{tabular}
 & $\longleftrightarrow$ &
\begin{tabular}{c}
\includegraphics{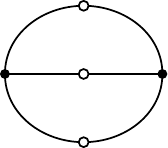}
\end{tabular}  \\ \\
\begin{tabular}{c}
\includegraphics{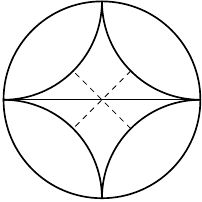} 
\end{tabular}
 & $\longleftrightarrow$ &
\begin{tabular}{c}
\includegraphics{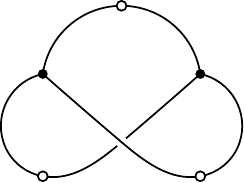}
\end{tabular}  \\ \\
\begin{tabular}{c}
\includegraphics{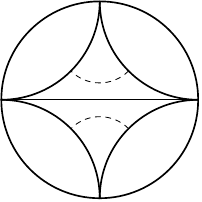} 
\end{tabular}
 & $\longleftrightarrow$ &
\begin{tabular}{c}
\includegraphics{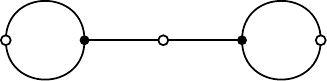}
\end{tabular}
\end{tabular}
\caption{The three triangular maps with two faces together with their incidence diagrams. The triangular map in the middle row corresponds to the only way to build a torus using two triangles, and the two others, top and bottom rows, correspond to the two different ways to build a sphere using two triangles.}
\label{tab:map:2}
\end{center}
\end{table}

\begin{table}
\begin{center}
\includegraphics{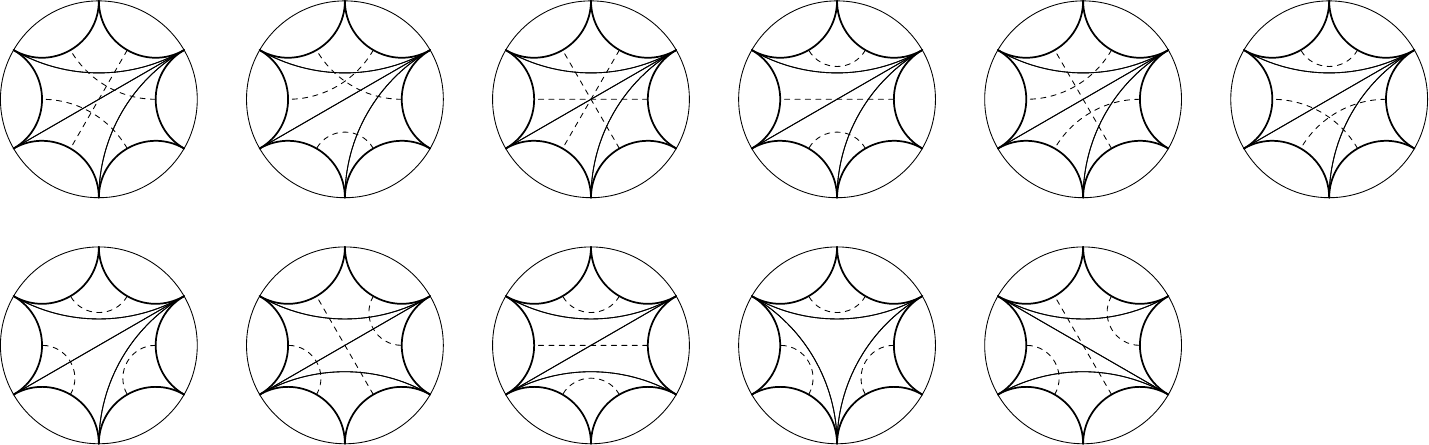}
\caption{The eleven triangular maps with four faces. According to table \ref{tab:unrooted:map:genus}, six of them are spheres while five of them are toruses. This is easily checked using Euler formula.}
\label{tab:map:4}
\end{center}
\end{table}

\begin{table}
\begin{center}
\includegraphics{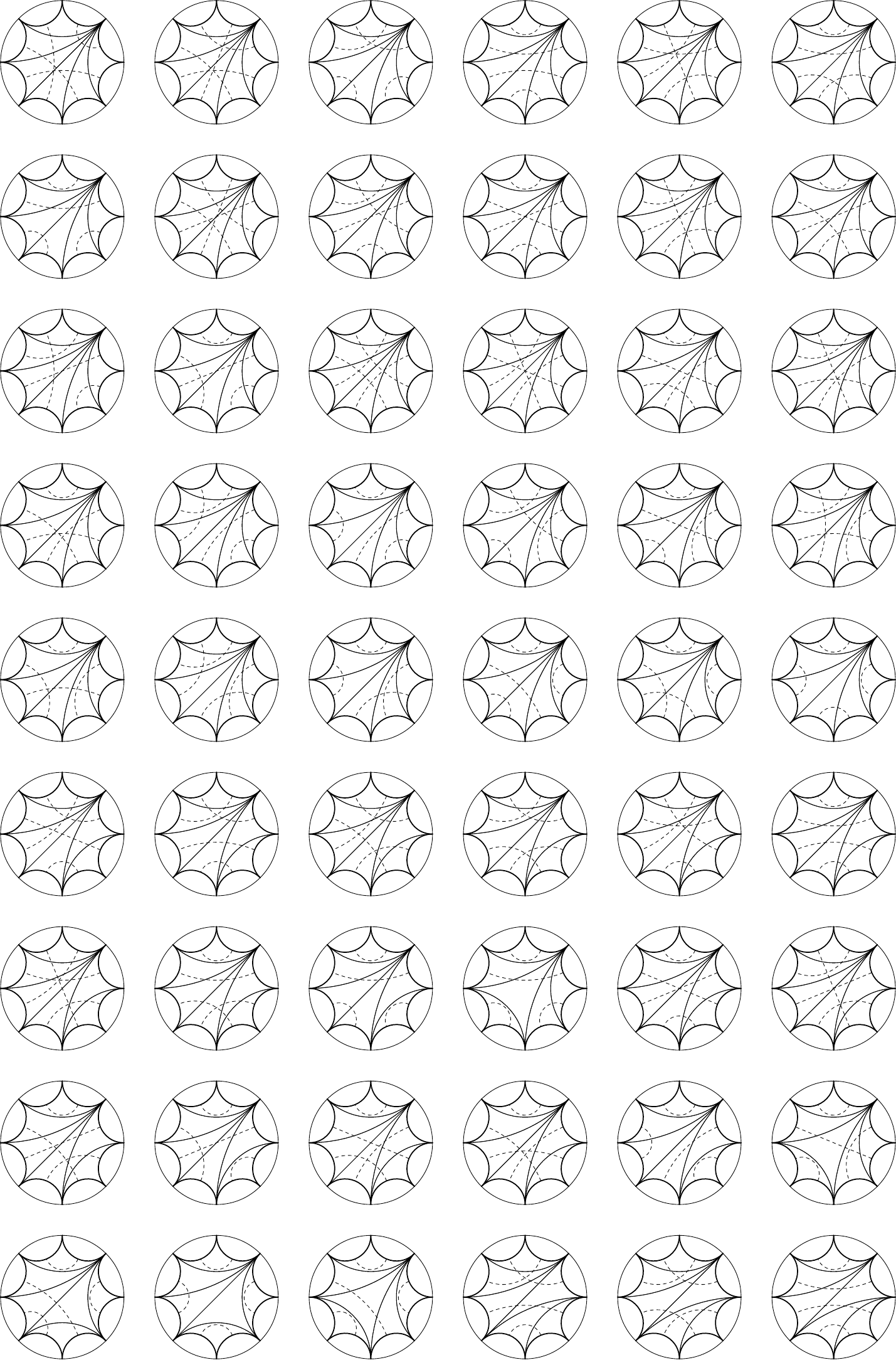}
\caption{The eighty one triangular maps with six faces (first part).}
\label{tab:map:6:1}
\end{center}
\end{table}

\begin{table}
\begin{center}
\includegraphics{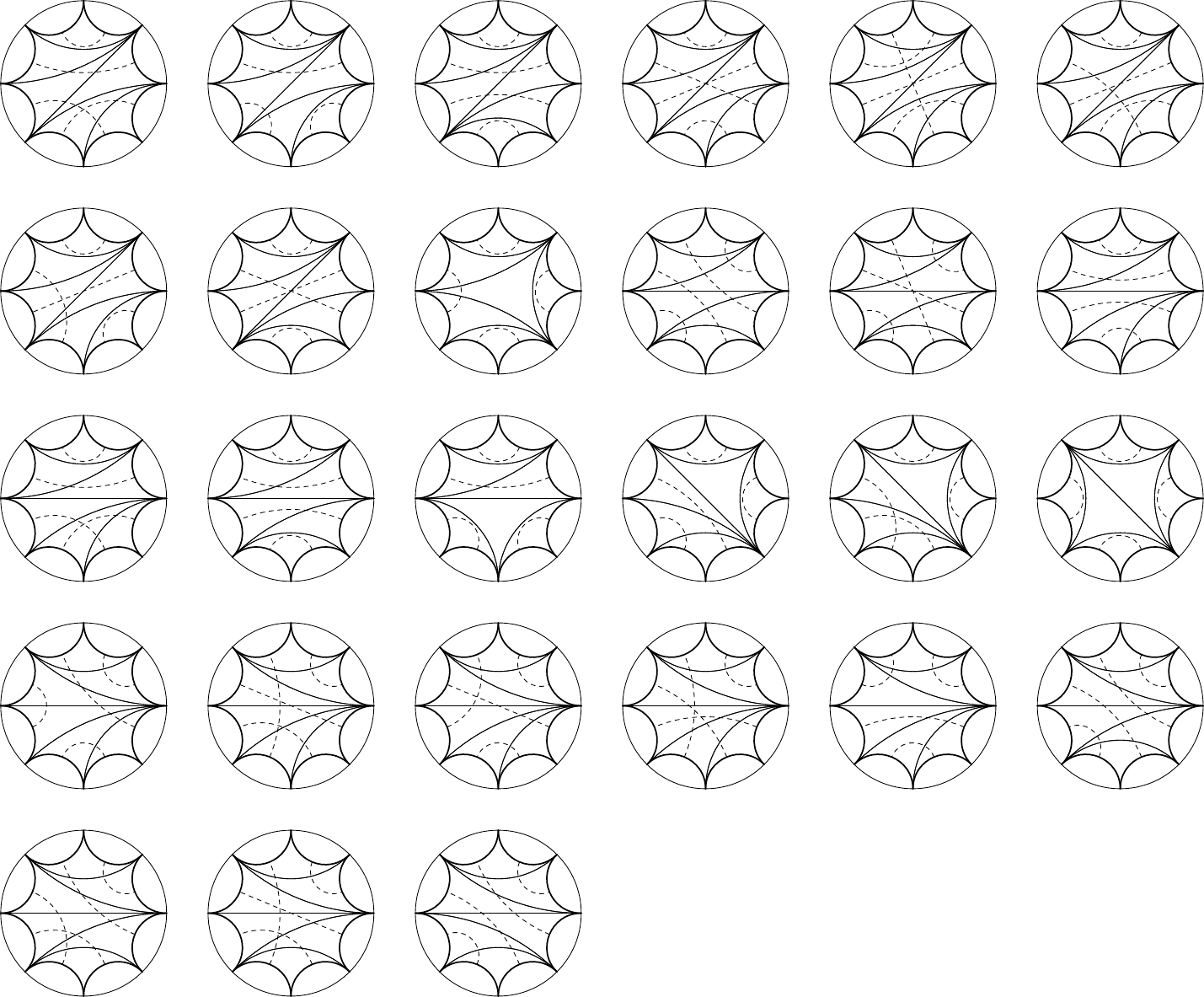}
\caption{The eighty one triangular maps with six faces (last part, continued from table \ref{tab:map:6:1}).}
\label{tab:map:6:2}
\end{center}
\end{table}


\begin{table}
\begin{center}
\begin{align*}
\tilde{M}^\bullet_3(t)\, &= 5\,{t}^{6}+60\,{t}^{12}+1105\,{t}^{18}+27120\,{t}^{24}+828250\,{t}^{30}+30220800\,{t}^{36}
\\&\quad+1282031525\,{t}^{42}+61999046400\,{t}^{48}+3366961243750\,{t}^{54}
\\&\quad+202903221120000\,{t}^{60}+13437880555850250\,{t}^{66}+970217083619328000\,{t}^{72}
\\&\quad+75849500508999712500\,{t}^{78}+6383483988812390400000\,{t}^{84}
\\&\quad+575440151532675686278125\,{t}^{90}+55318762960656722780160000\,{t}^{96}
\\&\quad+5649301494178851172304968750\,{t}^{102}
\\&\quad+610768380520654474629120000000\,{t}^{108}
\\&\quad+69692599846542054607811528918750\,{t}^{114}
\\&\quad+8370071726919812448859648819200000\,{t}^{120}+o(t^{120})
\end{align*}
\caption{Development of the generating series $\tilde{M}^\bullet_3(t)$, up to order one hundred and twenty. It gives as the coefficient of $t^{6n}$ the number of connected unrooted unlabeled triangular maps with $n$ arcs, thus $n/2$ undirected edges and $n/3$ triangular faces (A062980). If we denote by $a_n$ that coefficient, the recurrence is as follows: $a_1 = 5$ and  for $n \ge 1$, $a_{n+1}= (6n+6)\,a_n + \sum_{k = 1}^{n-1} a_k\,a_{n-k}$.}
\label{tab:ser:gen;t3dot:unlabeled}
\end{center}
\end{table}

\begin{table}
\begin{center}
\begin{align*}
\tilde{M}_3(t)\, &= 3\,{t}^{6}+11\,{t}^{12}+81\,{t}^{18}+1228\,{t}^{24}+28174\,{t}^{30}+
843186\,{t}^{36}+30551755\,{t}^{42}
\\&\quad+1291861997\,{t}^{48}+62352938720\,
{t}^{54}+3381736322813\,{t}^{60}
\\&\quad+203604398647922\,{t}^{66}+13475238697911184\,{t}^{72}+972429507963453210\,{t}^{78}
\\&\quad+75993857157285258473\,{t}^{84}+6393779463050776636807\,{t}^{90}
\\&\quad+576237114190853665462712\,{t}^{96}+55385308766655472416299110\,{t}^{102}
\\&\quad+5655262782600929403228668176\,{t}^{108}
\\&\quad+611338595145132827847686253456\,{t}^{114}
\\&\quad+69750597724332100283681465962492\,{t}^{120}+o(t^{120})
\end{align*}
\caption{Order one hundred and twenty development of the generating series $\tilde{M}_3(t)$ giving the number of connected unrooted unlabeled triangular maps with $n$ arcs, thus $n/2$ undirected edges and $n/3$ triangular faces (A129114).}
\label{tab:ser:gen:t3:unlabeled}
\end{center}
\end{table}

To adapt the generator algorithm to produce only regular rooted trivalent diagrams and thus rooted triangular maps, it suffice to remove the call to algorithms \ref{algo:tryclosedblack} and \ref{algo:tryclosedwhite} from line $6$ and $7$ of the $\textsc{Dispatch}$ procedure (algorithm \ref{algo:dispatch}) and the lines 2 to 5 from algorithm \ref{algo:generate}. Those two removals preserve the CAT property we thus get this way a constant amortized time generator for regular rooted triangular maps, as announced.

The tables \ref{tab:map:2} to \ref{tab:map:6:2} shows exhaustive lists of \emph{unrooted} triangular maps produced from the output of the generator algorithm. Those tables arise from a two steeps process. First, generate the full list of rooted regular diagrams of a given size, then, remove the duplicated diagrams from the list one obtain by forgetting the root. An easy and efficient way to do that is described in \cite{walsh83}, one has to put a linear order on the set of rooted trivalent diagrams of a given size, then one has to remove from the list the trivalent diagrams that are not minimal in their conjugacy class (two rooted diagrams being conjugated by definition, when they differ only by the position of their root).

The interpretation of the drawings of tables \ref{tab:map:2} to \ref{tab:map:6:2} deserves some explanation. For that purpose we adopt a geographical terminology. Ignoring for a moment the surrounding circle and the dashed lines of the drawings and the triangular regions are called the countries of the maps and the plain lines are the boundaries of their adjacent countries. One can distinguish the boundaries that are bordering two  
distinct countries (the inner boundaries) from those that are bordering a single country (the outer boundaries). The roads of the maps are symbolized by dashed lines connecting, in a two-by-two fashion, the outer boundaries of the maps.

To produce them, we have considered the planar rooted binary tree of the depth-first traversal of algorithm \ref{algo:visit:relabel}. As already noted, this algorithm provides natural cuttings for the associated trivalent diagrams. Those cuttings arise as what we previously called backward connections. In contrast, the edges of the traversal tree correspond to what we previously called forward connections.
In the graphical representation, we use an embedding of the produced triangulated polygon in the Poincar disc model of the hyperbolic plane as it seems the natural setting for generic non-overlapping triangular tilings. The surrounding circle around each figure is of course irrelevant to the structures.

\begin{table}
\begin{center}
\begin{tabular}{c rc | c r r r r r r r r c }
\hline \hline
 &&&&  2   &      4       &  6   &      8     &   10   &     12 &       14 & \\
\hline
 &0 &&&  4    &    32 &      336   &   4096  &  54912 &   786432  & 11824384 \\
& 1 &&&  1    &    28  &     664 &    14912  &  326496  & 7048192 & 150820608 \\
& 2 &&&     0    &     0     &  105   &   8112  &  396792 & 15663360 & 544475232 \\
 &3 &&&     0    &     0    &     0      &   0  &   50050  & 6722816 & 518329776 \\
 &4 &&&      0    &     0     &    0     &    0   &      0      &   0  & 56581525 \\
\hline
\hline
\end{tabular}
\caption{The number of \emph{rooted} triangular maps by genus (horizontally) and number of faces (vertically).}
\label{tab:rooted:map:genus}
\end{center}
\end{table}

\begin{table}
\begin{center}
\begin{tabular}{crc | c r r r r r r r rc }
\hline \hline
 &&&&  2   &      4       &  6   &      8     &   10   &     12 &       14& \\
\hline
& 0 &&& 2   &      6   &    26   & 191  &   1904  &      22078 &       282388    \\
& 1 &&& 1   &      5   &    46   & 669  & 11096  &   196888  &    3596104    \\
& 2 &&& 0   &      0   &      9  &  368  &  13448  &  436640  &  12974156    \\
& 3 &&& 0   &      0   &      0   &      0  &    1726  &  187580  &  12350102    \\
& 4 &&& 0   &      0   &      0   &      0  &           0  &             0   &    1349005    \\
\hline
\hline
\end{tabular}
\caption{The number of \emph{unrooted} triangular maps by genus (horizontally) and number of faces (vertically).}
\label{tab:unrooted:map:genus}
\end{center}
\end{table}

The tables \ref{tab:ser:gen;t3dot:unlabeled} and \ref{tab:ser:gen:t3:unlabeled} give the number of rooted triangular maps and unrooted triangular maps in the form of generating series. They are also part of \cite{sloan:oeis} under the references (A062980) and (A129114). Their computation is very similar to that of tables \ref{tab:ser:gen:d3dot:unlabeled} and \ref{tab:serie:diag:triv}, which is explained in detail in \cite{vidal06}. We shall explain in \cite{vidal07b} in a unified fashion how one can compute generating series for rooted and unrooted unlabeled maps of various kind and in \cite{vidal07c} the unexpected relation of this sequence to the asymptotcis of the Airy function.

As another byproduct of the exhaustive list obtained from the generating algorithm, one can get the precise number of rooted and unrooted triangular map having a given genus and a given number of triangular faces. Tables \ref{tab:rooted:map:genus} and \ref{tab:unrooted:map:genus} summarize those results for small number of faces.
Recently, M. Krikun \cite{krikun07} kindly communicated us recurrence relations satisfied by the entries of table \ref{tab:rooted:map:genus} which he obtained by a clever recursive decomposition of rooted triangular maps. Those recurrence relations make it  possible to evaluate easily those numbers without running the generator algorithm. Unfortunately, no general recurrence relation is known for the entries of table \ref{tab:unrooted:map:genus}. The thesis  \cite{walsh71} and articles \cite{walshlehan72, walshlehan75} contain the first enumerations of rooted maps of a given genus.

Some lines of those two tables were previously known. For instance, the first line of table \ref{tab:rooted:map:genus} is the number of spherical rooted triangular maps by the number of its faces \cite{mullin70}. The first line of table \ref{tab:unrooted:map:genus} is its unrooted counterpart. It is computed by impressive closed formulae in a recent paper by Liskovets, Gao and Wormald \cite{liskovets05}.
The diagonal terms of those two tables also received close attention. For instance, in \cite{zagier86} Harer and Zagier computed the Euler-Wall characteristic of the mapping class group of once pointed genus $g$ closed oriented surfaces by a remarkable combinatorial reduction of the problem in which rooted combinatorial maps with one vertex are counted by genus yielding the diagonal sequence of the first table 1, 105, 50050, 56581525, .... The diagonal sequence of the second table: 1, 9, 172, 1349005, .... gives the number of unrooted triangular maps of genus $g$ with only one vertex. It has been studied at depth in the article \cite{vdovina02} by A. Vdovina and R. Bacher.

\section{Concluding Remarks and Prespectives}

The generating algorithm presented in this paper (section \ref{sect:generator}) may lead to trivial adaptations to generate wider classes of diagrams and combinatorial maps, possibly with prescribed degree lists for vertices or faces. Basically, it can be simply generalized to produce any connected pair of permutations with prescribed cyclic types, up to simultaneous conjugacy.

Another way to extend the study, would be to modify the $\textsc{Dispatch}$ procedure (algorithm \ref{algo:dispatch} of section \ref{sect:generator}) to generate not an exhaustive cover of the partial cases, but instead a single case of them picked at random. This would result, in a fairly straightforward fashion, in a random sampler algorithm of the corresponding combinatorial structures instead of an exhaustive generator. The difficulty there, is to precompute precise conditional probability tables in order to control the probability distribution of the generated structures by bayesian techniques.
Such tables of conditional probabilities could be computed with the help of generating series techniques namely by following the particular recursive structure of the algorithm and translating this recursive structure in functional equations on the generating series. This appeal for a further investigation and could be dealt with in a subsequent paper.

If one designs the Dispatch procedure to pick at random a case with uniform probability this would unfortunately not result in a uniform distribution of the output structures but it can still prove useful to produce test cases to many algorithms operating on triangulations. A much beter approche but still not perfect is to generate at random the two permutations $\sigma_\bullet$ and $\sigma_\circ$ having the right cycle types with uniform distribution and rejecting each time the resulting map is not connected. This procedure produces each map with probability proportional to $1/A$ where $A$ is the number of its automorphisms. Since the maps have no automorphism other than identity in the vast majority of cases, that statistical bias is in practice unnoticeable.

\section{Acknowledgements}

I am grateful to professors D. Bar Nathan, P. Flajolet, F. Hivert, M. Huttner, M. Krikun, M. Petitot, B. Salvy, G. Shaeffer, N. Thiery, and D. Zvonkine for useful discussions and warm encouragements. Comments and suggestions from the anonymous referees helped a lot to improve the clarity of the paper and the quality of presentation.

\nocite{*}
\bibliography{art}
\bibliographystyle{plain}

\end{document}